\def\Extended{extended}
\documentclass[graybox]{svmult}

\usepackage{type1cm}        
\usepackage{makeidx}         
\usepackage{graphicx}        
\usepackage{multicol}        
\usepackage[bottom]{footmisc}

\usepackage{newtxtext}       %
\usepackage[varvw]{newtxmath}       

\makeindex             

\newif\ifextended
\ifx\Extended\Kaba\else\extendedtrue\fi
\newif\ifprivate
\ifx\Private\Kaba\else\privatetrue\extendedtrue\fi
\newif\iftalk
\ifx\Talk\Kaba\else\talktrue\fi
\ifextended
\usepackage[bookmarks=true,backref=page,bookmarksnumbered=true, 
  pdftitle={Reflection and Recurrence},%
  pdfauthor={Sakae Fuchino)},%
  pdfsubject={},%
  pdfkeywords={Generic large cardinal, Laver-generic large cardinal, Reflection Principles, Maximality Principles, Recurrence Axioms}]{hyperref}
\hypersetup{
  colorlinks=true,
  linkcolor=blue,
  filecolor=blue,      
  citecolor=cyan,
  urlcolor=cyan}
\usepackage{marginnote}
\usepackage[normalem]{ulem}
\newcommand\redout{\bgroup\markoverwith
{\textcolor{red}{\rule[.5ex]{2pt}{0.8pt}}}\ULon} 
\usepackage[T1]{fontenc} 
\fi
\newcount\minute	
\newcount\hour		
\newcount\hourMins  
\def\now%
{
  \minute=\time    
  \hour=\time \divide \hour by 60 
  \hourMins=\hour \multiply\hourMins by 60
  \advance\minute by -\hourMins 
  \zeroPadTwo{\the\hour}:\zeroPadTwo{\the\minute}(GMT)%
}
\def\zeroPadTwo#1%
{
  \ifnum #1<10 0\fi    
  #1
}
\def\daimaru#1{\makebox[1em][c]{\mbox{\leavevmode\lower.144ex\hbox{
        \rlap{\hbox to 
          0.76em{\hfil\mbox{}\hfill{}\raisebox{0.054ex}{\scalebox{1.2}{\begin{CJK}{UTF8}{goth}○\end{CJK}}}\hfil}}
        \raise0.342ex\hbox to 1em{\hfil{\hspace{0.16em}\footnotesize#1}\hfil}}}}\,}
\ifprivate
\newcommand{\Label}[1]{\label{#1}\marginnote{{\color{cyan}\renewcommand{\baselinestretch}{0.4}\tiny{#1}}}}

\newcommand{\Bibitem}[1]{\bibitem{#1}\marginnote{\color{cyan}%
    {\renewcommand{\baselinestretch}{0.4}\tiny 
		  {#1}}}}
\else
\newcommand{\Label}[1]{\label{#1}}

\newcommand{\Bibitem}[1]{\bibitem{#1}}
\fi
\newcommand{\pageof}[1]{p.\pageref{#1}}
\ifextended
\newcommand{\jumpto}[2]{\hyperlink{#1}{{\color{blue}#2\,}}}
\newcommand{\Hypertarget}[1]{\hypertarget{#1}}
\else
\newcommand{\jumpto}[2]{{#2}}
\newcommand{\Hypertarget}[1]{}
\fi
\def\memo#1{\ifprivate\marginnote{{\darkred\normalsize%
      \renewcommand{\baselinestretch}{0.4}\tiny\mbox{}\vspace{-2.52ex}
      \par\relax%
			#1\par\mbox{}}}\else\fi}%
\def\memoL#1{\ifprivate\marginnote{{\darkred\normalsize%
      \renewcommand{\baselinestretch}{0.4}\tiny\mbox{}\vspace{-2.52ex}
      \par\relax%
			#1\hfill\mbox{}}}\else\fi}%
\def\memox#1{}
\def\memoRx#1{}
\def\memoLx#1{}
\def\imemox#1{}

\usepackage{bbm} 
\usepackage{accents} 
\ifextended
\usepackage{CJKutf8} 
\usepackage[many]{tcolorbox} 
\fi
\usepackage{mathtools}
\newtagform{stareq}{($^*$}{)}
\newcommand{\utildeT}[1]{%
  \underaccent{\scalebox{0.98}{$\sim$}}{#1}}
\newcommand{\utildeS}[1]{%
	\hbox to 0pt{\smash{$\mathop{\scriptstyle #1}\limits_{%
				\raisebox{0.6ex}[0pt]{$\scriptstyle\sim$}}$}\hss}%
	\relax\phantom{\mathord{{#1}_{\rule[-0.6ex]{0pt}{1pt}}}}}
\newcommand{\utildeSS}[1]{%
	\hbox to 0pt{$\mathop{\scriptscriptstyle #1}%
		\limits_{\scriptscriptstyle\sim}$\hss}%
		\relax\phantom{\underline{#1}}}
\newcommand{\utilde}[1]{%
	\mathchoice{\utildeT{#1}}{\utildeT{#1}}{\utildeS{#1}}{\utildeSS{#1}}}
\newcommand{\bbd}[1]{{\mathbbm{#1}}}
\newcommand{\setof}[2]{\{#1\,:\,#2\}}
\newcommand{\ssetof}[1]{\{#1\}}
\newcommand{\seqof}[2]{\langle#1\,:\,#2\rangle}
\newcommand{\pairof}[1]{\langle#1\rangle}
\newcommand{\psof}[1]{{\mathcal P}\/(#1)}
\newcommand{\Pkl}[2]{\ifx\bakakaba#1\bakakaba\ifx\bakakaba#2\bakakaba{\mathcal 
    P}_\kappa(\lambda)\else{\mathcal P}_\kappa(#2)\fi\else{\mathcal P}_{#1}(#2)\fi}
\newcommand{\cardof}[1]{\mathopen{|\,}#1\mathclose{\,|}}
\newcommand{\forces}[2]{\,\|\hspace{-.35ex}\mbox{\sf--}_{\,#1\,}%
\mbox{\rm``}\,#2\,\mbox{\rm''}}

\newcommand{\modelof}[1]{\models\!\mbox{\rm``\,}#1\mbox{\rm''}}
\newcommand{\circleq}{\mathrel{{\leqslant}%
		\hspace{-0.8ex}{\lower-0.42ex\hbox{$\scriptscriptstyle\circ$}}}}

\newcommand{\unmeralof}[1]{\llcorner\hspace{-0.02em}#1\hspace{-0.02em}\lrcorner}
\newcommand{\mapping}[3]{#1:#2\rightarrow #3}
\newcommand{\fnsp}[2]{\mbox{}^{{#1}\hspace{-0.02em}}#2}
\newcommand{\imageof}{{}^{\,{\prime}{\prime}}}

\newcommand{\Elembed}[4]{#1:#2\stackrel{\raisebox{-0.8ex}{\smash{\scalebox{0.7}{$\prec$}}\hspace{0.8ex}}}{\rightarrow}_{#4}#3}
\newcommand{\checked}{\hspace{0.02em}^{\checkmark}}
\newcommand{\continuum}{2^{\aleph_0}}

\newcommand{\bbone}{{\mathord{\mathbbm{1}}}}

\newcommand{\ol}[1]{\overline{#1}}
\newcommand{\restr}{\restriction}
\newcommand{\trcl}{\mathop{\mbox{\it trcl\/}}}

\newcommand{\ctenten}{,\mbox{}\hspace{0.08ex}{.}{.}{.}\hspace{0.1ex}}
\newcommand{\ctentenc}{,{}\linebreak[0]\hspace{0.04ex}{{.}{.}{.}\hspace{0.1ex},\,}%
  \linebreak[0]}

\newcommand{\xmbox}[1]{ $\relax{\rm #1}\relax$ }

\newcommand{\BaB}{\bbd{B}}

\newcommand{\poP}{\bbd{P}}
\newcommand{\poQ}{\bbd{Q}}

\newcommand{\poS}{\bbd{S}}
\newcommand{\uta}{\utilde{a}}
\newcommand{\utb}{\utilde{b}}
\newcommand{\utpoP}{\utilde{\mathbbm{P}}}
\newcommand{\utpoQ}{\utilde{\mathbbm{Q}}}
\newcommand{\utpoR}{\utilde{\mathbbm{R}}}
\newcommand{\genF}{\mathbbm{F}}
\newcommand{\genG}{\mathbbm{G}}
\newcommand{\genH}{\mathbbm{H}}

\newcommand{\condp}{\mathbbm{p}}

\newcommand{\uniU}{{\sf U}}
\newcommand{\uniV}{{\sf V}}
\newcommand{\uniW}{{\sf W}}
\newcommand{\gmA}{\mathfrak{A}}
\newcommand{\gmB}{\mathfrak{B}}

\newcommand{\Lin}{{\calL}_{\in}}
\newcommand{\ZF}{{\sf ZF}}
\newcommand{\ZFC}{{\sf ZFC}}
\newcommand{\CH}{{\sf CH}}
\newcommand{\MA}{{\sf MA}}
\newcommand{\MM}{{\sf MM}}
\newcommand{\SCH}{{\sf SCH}}
\newcommand{\PFA}{{\sf PFA}}
\newcommand{\DRP}{{\sf DRP}}
\newcommand{\GRP}{{\sf GRP}}
\newcommand{\FRP}{{\sf FRP}}
\newcommand{\IC}{{\sf IC}}
\newcommand{\IU}{{\sf IU}}
\newcommand{\RP}{{\sf RP}}
\newcommand{\RcA}{{\sf RcA}}
\newcommand{\RcAp}{{\sf RcA$^+$}}
\newcommand{\IMH}{{\sf IMH}}
\newcommand{\IGH}{{\sf IGH}}
\newcommand{\LGM}{{\sf LGM}}
\newcommand{\MP}{{\sf MP}}
\newcommand{\On}{{\sf On}}
\newcommand{\Card}{{\sf Card}}

\newcommand{\crit}{\mbox{\it crit\/}}

\newtheorem{Lemma}[theorem]{Lemma}
\newtheorem{Prop}[theorem]{Proposition}
\newtheorem{Cor}[theorem]{Corollary}
\newtheorem{Claim}{Claim}[theorem]
\newtheorem{ThmA}{Theorem A\!\!}
\newtheorem{CorA}[ThmA]{Corollary A\!\!}
\ifextended
\tcolorboxenvironment{theorem}{
  colback=blue!4!white,
  boxrule=0pt,
  boxsep=0pt,
  left=6pt,right=6pt,top=2pt,bottom=4pt,
  oversize=4pt,
  sharp corners,
  before skip=\topsep,
  after skip=\topsep,
  breakable
}
\tcolorboxenvironment{Prop}{
  colback=blue!4!white,
  boxrule=0pt,
  boxsep=0pt,
  left=6pt,right=6pt,top=2pt,bottom=4pt,
  oversize=4pt,
  sharp corners,
  before skip=\topsep,
  after skip=\topsep,
  breakable
}
\tcolorboxenvironment{Lemma}{
  colback=blue!4!white,
  boxrule=0pt,
  boxsep=0pt,
  left=6pt,right=6pt,top=2pt,bottom=4pt,
  oversize=4pt,
  sharp corners,
  before skip=\topsep,
  after skip=\topsep,
  breakable
}
\tcolorboxenvironment{Cor}{
  colback=blue!4!white,
  boxrule=0pt,
  boxsep=0pt,
  left=6pt,right=6pt,top=2pt,bottom=4pt,
  oversize=4pt,
  sharp corners,
  before skip=\topsep,
  after skip=\topsep,
  breakable
}
\tcolorboxenvironment{CorA}{
  colback=blue!4!white,
  boxrule=0pt,
  boxsep=0pt,
  left=6pt,right=6pt,top=2pt,bottom=4pt,
  oversize=4pt,
  sharp corners,
  before skip=\topsep,
  after skip=\topsep,
  breakable
}
\tcolorboxenvironment{Claim}{
  colback=pink!5!white,
  boxrule=0pt,
  boxsep=0pt,
  left=6pt,right=6pt,top=2pt,bottom=4pt,
  oversize=4pt,
  sharp corners,
  before skip=\topsep,
  after skip=\topsep,
  breakable
}
\tcolorboxenvironment{itemize}{
  colback=yellow!5!white,
  boxrule=0pt,
  boxsep=0pt,
  left=1pt,right=1pt,top=2pt,bottom=2pt,
  oversize=2pt,
  sharp corners,
  before skip=\topsep,
  after skip=\topsep,
  breakable
}

\fi
\newcommand{\ubecause}[3]{\underbrace{{}#1{}%
  \ifx\bakakaba#2\bakakaba\rule[-0.72ex]{0pt}{1pt}\else\rule[#2]{0pt}{1pt}\fi}_{\mbox{\footnotesize\clap{#3}}}}
\newcommand{\obecause}[3]{\overbrace{{}#1{}%
  \ifx\bakakaba#2\bakakaba\rule[1.62ex]{0pt}{1pt}\else\rule[#2]{0pt}{1pt}\fi}^{\mbox{\footnotesize\clap{#3}}}}
\newcommand{\prf}{\noindent{\it Proof.\ \ }\ignorespaces}
\newcommand{\prfofClaim}{\raisebox{-.4ex}{\Large $\vdash$\ \ }}
\renewcommand{\qed}{\mbox{}\hfill{\large$\Box$}}
\newcommand{\qedofClaim}%
{\mbox{}\hfill\raisebox{-.4ex}{\Large $\dashv$ }\nolinebreak%
\mbox{\tiny~(Claim~\number\theClaim)}}
\newcommand{\qedof}[1]%
           {\mbox{} \hspace*{\fill}{{\Large$\Box$}{\tiny~(#1)}}}
\newcommand{\Qedof}[1]%
{\mbox{} \hspace*{\fill}{{\large$\Box$}%
{\tiny~(#1~\number\thetheorem)}}}
\newcommand{\QedofL}[1]%
{\mbox{} \hspace*{\fill}{{\large$\Box$}%
{\tiny~(#1~\number\thelemma)}}}
\newcommand{\QedofTA}[1]%
{\mbox{} \hspace*{\fill}{{\large$\Box$}%
{\tiny~(#1~\number\theThmA)}}}
\newcommand{\theoremof}[1]{Theorem \ref{#1}}
\newcommand{\lemmaof}[1]{Lemma \ref{#1}}
\newcommand{\Lemmaof}[1]{Lemma \ref{#1}}
\newcommand{\Propof}[1]{Proposition \ref{#1}}
\newcommand{\Corof}[1]{Corollary \ref{#1}}
\newcommand{\Claimof}[1]{Claim \ref{#1}}
\newcommand{\sectionof}[1]{Section \ref{#1}}
\newcommand{\footnoteof}[1]{Footnote \ref{#1}}
\newcommand{\qedoftheorem}{\Qedof{Theorem}}
\newcommand{\qedofLemma}{\Qedof{Lemma}}
\newcommand{\qedofProp}{\Qedof{Proposition}}

\newcommand{\qedofCor}{\Qedof{Corollary}}
\newcommand{\qedskip}{\medskip}
\newcommand{\qedofCorA}{\QedofTA{Corollary A\!\!}}

\newenvironment{xitemize}{\begin{list}{}{\parsep=0.5\smallskipamount%
			\itemindent=-0.1ex%
			\itemsep=0.5\smallskipamount\leftmargin=3.0em\labelwidth=3em\labelsep=0.5em}}%
							 {\end{list}}
\newcommand{\equationof}[1]{($^*$\ref{#1})}
\def\xitem[#1]{\refstepcounter{equation}{%
  \def\baka{#1}\ifx\baka\empty\else\label{#1}\fi}%
  \item[{\rm\makebox[3.2em][c]{($^*$\arabic{equation})}}]
   \ifprivate\marginnote{\mbox{}\hfill\darkred\tiny #1}\ignorespaces\fi}
\def\xitemas[#1]{\item[{\rm\makebox[3.2em][c]{($^*$\ref{#1})}}]
   \ifprivate\marginnote{\mbox{}\hfill\darkred\tiny #1}\ignorespaces\fi}
\newcounter{frmla}
\def\xitemA[#1]{\refstepcounter{frmla}{%
  \def\baka{#1}\ifx\baka\empty\else\label{#1}\fi}%
  \item[{\rm\makebox[3.2em][c]{($^\aleph$\arabic{frmla})}}]
   \ifprivate\marginnote{\mbox{}\hfill\darkred\tiny #1}\ignorespaces\fi}
\newcommand{\xitemof}[1]{($^*$\ref{#1})}
\def\xitemAof#1{{\rm({$^\aleph$\ref{#1}})}}

\newcommand{\calD}{{\mathcal D}}

\newcommand{\calF}{{\mathcal F}}

\newcommand{\calH}{{\mathcal H}}

\newcommand{\calL}{{\mathcal L}}

\newcommand{\calP}{{\mathcal P}}
\newcommand{\calQ}{{\mathcal Q}}

\newcommand{\calS}{{\mathcal S}}

\newcommand{\Col}{{\rm Col}}
\newcommand{\mahlo}{{\it M\ell\,}}

\newcommand{\natnums}{{\bbd{N}}}
\newcommand{\LT}{{<}\,}
\newcommand{\LE}{{\leq}\,}

\ifextended\newcommand{\darkred}{\color[rgb]{0.8,0.1,0.1}}\else\newcommand{\darkred}{}\fi
\ifextended\newcommand{\extendedcolor}{\color[rgb]{0.33, 0.41, 0.47}}\else
\newcommand{\extendedcolor}{}\fi
\ifextended\newcommand{\It}{\darkred\it}\else\newcommand{\It}{\it}\fi
\def\assert#1{\noindent\makebox[4.8ex][r]{\rm(\makebox[2.2ex][c]{#1})}\ \ \ignorespaces}
\def\wassert#1{\assert{#1}}
\def\wassertof#1{\makebox[4.8ex][r]{\rm(\makebox[2.2ex][c]{#1})\ }}%
\def\assertof#1{{\rm(#1)}}%
\def\glqq{\mbox{,\hspace{-.1em},}}
\def\grqq{\mbox{`\hspace{-.05em}`}}
\newcommand{\st}{such that}
\newcommand{\wrt}{with respect to}
\newcommand{\Wolog}{Without loss of generality}

\newcommand{\po}{poset}
\newcommand{\pos}{posets}
\newcommand{\refl}{{\mathfrak{r}\mathfrak{e}\mathfrak{f}\mathfrak{l}\,}}

\begin{document}

\title*{Reflection and Recurrence.}
\ifextended
\author{Sakaé Fuchino
}
\else
\author{Sakaé Fuchino
}
}
\fi

\institute{Sakaé Fuchino \at Graduate School of System Informatics, Kobe University, 
  Nada, Kobe 657-8501 Japan,\\\email{fuchino@diamond.kobe-u.ac.jp}}
%
%
\maketitle

{\def\thefootnote{}
\ifextended \footnotetext{\hspace*{-0.3em}\indent This is an extended 
  postprint version of  
  the paper with the same  
  title which will appear in a Festschrift in occasion of the 75.\ birthday of Professor 
  Janos Makowsky. The most recent file of this version is downloadable as: \\
  \scalebox{0.94}[1]{\url{https://fuchino.ddo.jp/papers/reflection\_and\_recurrence-Janos-Festschrift-x.pdf}}\smallskip\\
  \indent
  All additional texts not included in the published version and all major corrections 
  after the publication are 
  marked with\quad {\extendedcolor ``dark electric blue'' foreground color}. 
}
\else \footnotetext{\Label{postprint}\hspace*{-0.3em}An extended postprint version of the chapter is 
  downloadable as: \\ 
  \url{https://fuchino.ddo.jp/papers/reflection\_and\_recurrence-Janos-Festschrift-x.pdf}}
\fi}
\ifextended
\phantomsection
\addcontentsline{toc}{section}{Abstract}
\fi
\abstract{ We examine the Zermelo Fraenkel set theory with Choice (\ZFC) enhanced by one of 
  the (structural) reflection principles down to a small 
  cardinal and/or Recurrence Axioms defined below. The strongest forms of 
  reflection principles spotlight the three scenarios in which 
  the size of the continuum is 
  either $\aleph_1$, or $\aleph_2$, or very large, while the maximal setting of Recurrence 
  Axioms points to the set-theoretic universe with the continuum of size $\aleph_2$.
  \newline\indent
  We discuss that both the Reflection Principles and Recurrence Axioms can be 
  construed as preferable candidates of the extension of \ZFC\ in terms of the criteria 
  of Gödel's Program. From this view point, the maximal possible 
  (consistent)  combination of these principles and  
  axioms, or even some natural strengthening of the combination (which we want to 
  call ``Laver-generic Maximum'' (\LGM)) 
  may be considered as the ultimate extension of\ 
  \ZFC\ (of course ``ultimate'' only for now --- because of the Incompleteness Theorems): 
  \LGM\ resolves the size of the continuum to be $\aleph_2$ and integrates 
  practically all known statements consistent with \ZFC\ in itself either as its 
  consequences (as it is the case with Martin's Maximum$^{++}$) or as theorems holding in 
  many grounds of the universe (as it is the case with Cichoń's Maximum). 
  {\ifextended\small\extendedcolor\bigskip\newline \begin{CJK}{UTF8}{ipxm}{\bf 
        アブストラクト}\quad 以下で，
    選択公理付きの Zermelo  
    Fraenkel 集合論\ (\ZFC) を，小さな基数への (構造に関する) 反映原理 (Reflection 
    Principles) たちのうちの一つ，かつ／または， 
    回帰公理 (Recurrence Axioms) たちの一つで拡張した体系について考察する．反映原理のうちの
    一番強いものたちは，連続体のサイズが $\aleph_1$ であるか，$\aleph_2$ であるか，非常に
    大きくなるかという３つのシナリオにスポットライトをあてるのに対し，Recurrence Axioms の
    一番強いものたちの組合せは，
    連続体のサイズが $\aleph_2$ になる集合論の宇宙を指ししめしているようである．\
    \newline\indent 反映原理も回帰公理も，ゲーデルのプログラムの意味で，\ZFC\ の望ましい
    拡張を与えるものの候補と看倣せるので，これらの原理と公理の無矛盾なもののうちの
    極大な組合せは，\ZFC\ の究極の拡張と考えることができる --- これを我々は ``Laver-generic 
    Maximum'' (\LGM) と呼びたい．もちろんこれは「究極の」，とは言っても不完全性定理のために，
    現在のところの究極でしかないものなのだが，\LGM\ は，連続体問題を
    連続体濃度が $\aleph_2$ である，として解決し，現在までに知られている実質的に全ての\ 
    \ZFC\ と無矛盾な命題を，(たとえば \MM$^{++}$ がそうであるように) この拡張から導かれる
    定理として，あるいは，(たとえば Cichoń's Maximum がそうであるように) 多くの，
    集合論宇宙の\ 
    grounds での定理として，この理論に統合する．
    \end{CJK}
     \fi}}
{\iftrue
\noindent
{\ifextended\extendedcolor (\today\ \now\ version)\fi}
\phantomsection
\addcontentsline{toc}{section}{Contents}
\newcommand{\myscalebox}[1]{\scalebox{0.88}[1.06]{#1}}
\begin{quotation}
	\footnotesize
	\noindent
	\centerline{
      \normalsize\tt\quad\ Contents\hspace{6em}\mbox{}}\mbox{}\\
\ifextended\newcommand{\hhhyperref}[2]{\hyperref[#1]{#2}}\else
\newcommand{\hhhyperref}[2]{{#2}}\fi
{\mbox{}\hspace{0em}\tt\makebox[3.4ex][l]{\ref{intro}.}%
  \hhhyperref{intro}{\tt\myscalebox{Introduction}}}\ \ 
\dotfill\ \ {\pageref{intro}}\\   
{\mbox{}\hspace{0em}\tt\makebox[3.4ex][l]{\ref{sec:1}.}%
  \hhhyperref{sec:1}{\tt\myscalebox{Reflection down to a small cardinal}}}\ \ 
\dotfill\ \ {\pageref{sec:1}}\\   
{\mbox{}\hspace{0em}\tt\makebox[3.4ex][l]{\iftalk\hspace{-1em}▷\fi\ref{sec:2}.}%
  \hhhyperref{sec:2}{\tt\myscalebox{Recurrence, Maximality, and the solution(s) of the 
      Continuum Problem}}}\ \ \dotfill\ \ {\pageref{sec:2}}\\   
{\mbox{}\hspace{0em}\tt\makebox[3.4ex][l]{\ref{sec:2-0}.}%
  \hhhyperref{sec:2-0}{\tt\myscalebox{Restricted Recurrence Axioms}}}\ \ 
\dotfill\ \ {\pageref{sec:2-0}}\\    
{\mbox{}\hspace{0em}\tt\makebox[3.4ex][l]{\ref{sec:3}.}%
  \hhhyperref{sec:3}{\tt\myscalebox{Recurrence, Laver-generic large 
      cardinal, and beyond}}}\ \ \dotfill\ \ {\pageref{sec:3}}\\   
{\mbox{}\hspace{0em}\tt\makebox[3.4ex][l]{\ref{sec:4}.}%
  \hhhyperref{sec:4}{\tt\myscalebox{Toward the Laver-generic Maximum}}}\ \ 
\dotfill\ \ {\pageref{sec:4}}\\   
{\mbox{}\hspace{0em}\tt\makebox[3.4ex][l]{\ref{sec:5}.}%
  \hhhyperref{sec:5}{\tt\myscalebox{More about consistency strength}}}\ \ \dotfill\ \ {\pageref{sec:5}}\\  
{\mbox{}\hspace{0em}\hhhyperref{ref}{\tt References}}\ \ \dotfill\ \ 
  {\pageref{ref}}
\end{quotation}\bigskip\bigskip
\fi}

\usetagform{stareq} 
\section{Introduction}\Label{intro}
In the following, we examine the Zermelo Fraenkel set theory with the Axiom 
  of Choice (\ZFC) \memox{\normalsize!!!}
enhanced by one of 
the (structural) reflection principles down to a small 
cardinal and/or Recurrence Axioms defined below. The strongest forms of 
reflection principles (existence of a/the $\calP$-Laver-generic large cardinal --- see 
\sectionof{sec:1} below) spotlight 
the three scenarios in which the size of the continuum is 
either $\aleph_1$, or $\aleph_2$, or very large (see Theorems \ref{p-2}, \ref{p-3} and 
\ref{p-3-0}), while the maximal setting of Recurrence  
Axioms points to the set-theoretic universe with the continuum of size $\aleph_2$ (see the 
end of  
\sectionof{sec:2}).

As we are going to discuss in sections \ref{sec:1} and \ref{sec:2}, both of the Reflection 
Principles and Recurrence Axioms can be 
construed as preferable candidates of the extension of \ZFC\ in terms of 
Gödel's Program  (Gödel \cite{goedel}, see also Bagaria \cite{bagaria}). 
From this point of view the maximal possible 
(consistent) combination of these principles and 
axioms, or even some natural strengthening of the combination, that is, either the principle 
\LGM\ proposed in \sectionof{sec:4} or some further extension of it in the future 
(which we want to call ``Laver-generic Maximum'' (\LGM), see \pageof{LGM}) 
may be considered as the ultimate extension of\ 
\ZFC\ (of course ``ultimate'' only for now --- because of the Incompleteness Theorems): 
\LGM\ 
resolves the size of the continuum to be $\aleph_2$ and integrates 
practically all known statements consistent with \ZFC\ in itself either as its 
consequences (as it is the case with Martin's Maximum$^{++}$, see \theoremof{p-1}) or as theorems 
holding in many grounds of the universe  
(as it is the case with Cichoń's Maximum, Goldstern, Kellner and Shelah\cite{GKS}, 
Goldstern, Kellner, Mejía and Shelah \cite{GKMS}, see around 
\pageof{CichonM}) --- for more detailed discussions about the significance of these 
axioms,  see also the end of both of the Sections 
\ref{sec:1}, and \ref{sec:4}. 
\section{Reflection down to a small cardinal}
\Label{sec:1}
The small cardinal mentioned in the title of this section may be considered as not very 
small by non-set-theoretic mathematicians: It is known that many ``mathematical'' 
reflection statements with reflection number less than or equal to 
$\kappa_\refl:=\max\ssetof{\aleph_2,\continuum}$ hold (often in some extension of \ZFC). Some 
of them are even theorems in \ZFC. For example, 
\begin{theorem}\Label{p-a} \wassertof{1} {\rm(Dow \cite{dow})} If $X$ is a countably compact 
  Hausdorff non-metrizable space then there is a subspace $Y$ of $X$ of cardinality 
  $<\aleph_2$ \st\ $Y$ is also non-metrizable. \smallskip

  \wassert{2} Let $L(Q)$ be a logic with new (first-order) quantifier \st\ ``$Qx$ ...'' is 
  interpreted as ``there are uncountably many $x$ \st\ ...''. For any structure $\gmA$ of 
  countable signature, there is $\gmB\prec_{L(Q)}\gmA$ of size $<\aleph_2$. \qed
\end{theorem}

From very early on, it was known that, starting from a very large cardinal, we can 
construct models of set theory in which various strong statements on (structural) 
reflection down to $\LT\kappa_\refl$ hold. 

For example, Ben-David \cite{ben-david} in 1978 mentions a theorem by Shelah which states:
\begin{theorem}\Label{shelah-0}{\rm (S.\ Shelah, \cite{ben-david})} Suppose that $\kappa$ is supercompact and
  $\poP=\Col(\aleph_1,\kappa)$. Then, for $(\uniV,\poP)$-generic $\genG$, we have 
  \begin{equation*}
    \uniV[\genG]\models
    \begin{array}[t]{@{}l}
      \mbox{ for any structure }\gmA\mbox{ of countable signature, there 
        is }\gmB\prec_{L_{stat}}\!\gmA\\
      \mbox{of cardinality}<\aleph_2.
    \end{array}
  \end{equation*}
\end{theorem}
Here, $L_{stat}$ denotes the stationary logic with monadic second-order variables $X$ which 
run over countable subsets of the underlying set of respective structures and with the 
second-order quantifier $stat\,X$ which is to be interpreted as "there are stationarily 
many countable sets X". 

The elementary submodel relation $\gmB\prec_{L_{stat}}\!\gmA$ is defined by
$\gmB\models\varphi(b_0\ctenten)$\ \ $\Leftrightarrow$\ \ 
$\gmA\models\varphi(b_0\ctenten)$ for all $L_{stat}$-formula
$\varphi=\varphi(x\ctenten)$ without free second-order variables, and for all
$b_0\cdots\in|\gmB|$. 

Today, we can understand Shelah's theorem above as a special case of the following theorem.
For a class $\calP$ of \pos, a cardinal $\kappa$ is said to be {\It$\calP$-generically 
supercompact\/}  
if, for any $\lambda>\kappa$, there is a \po\ $\poP\in\calP$ \st, 
for $(\uniV,\poP)$-generic $\genG$, there are $j$, $M\subseteq\uniV[\genG]$ \st\
$\Elembed{j}{\uniV}{M}{\kappa}$,\footnote{With ``$\Elembed{j}{\uniV}{M}{\kappa}$'' we 
  denote the situation that $M$ is transitive, $j$ is an elementary embedding of $\uniV$ 
  into $M$, and $\kappa$ is the critical point of $j$.}, $j(\kappa)>\lambda$, and
$j\imageof{\lambda}\in M$. 

\begin{theorem}\Label{p-0}
  Suppose that $\aleph_2$ is $\calP$-generically supercompact where $\calP$ is the class of 
  all 
  $\LT\aleph_1$-closed \pos. Then
  \begin{xitemize}
  \xitem[eq:-0]
    for any structure $\gmA$ of countable signature, there 
    is $\gmB\prec_{L_{stat}}\!\gmA$
    of cardinality $\LT\aleph_2$.
  \end{xitemize}
\end{theorem}
\prf The condition ``$\aleph_2$ is $\calP$-generic supercompact for $\calP=$ the class of all
$\LT\aleph_1$-closed \pos.'' is equivalent to the Game Reflection Principle (\GRP)  
(Theorem 8 in König \cite{koenig} --- see Fuchino, Ottenbreit Maschio Rodrigues and Sakai 
\cite{I} for a generalization of König's  
theorem in \cite{koenig} --- note 
that what we call \GRP\ here and \cite{I} is called ``the global Game Reflection 
Principle'' in \cite{koenig}). By Theorem 4.7 in \cite{I}, \GRP\ implies \equationof{eq:-0}. 
\qedoftheorem\qedskip

The downward Löwenheim-Skolem Theorem \equationof{eq:-0} for $L_{stat}$ is actually a 
strong reflection property. For example the reflection of uncountable coloring number of 
graphs down to $\LT\aleph_2$ (the following \equationof{eq:-1}) is a consequence of 
\equationof{eq:-0}:

\begin{xitemize}
\xitem[eq:-1]
  For any graph $G$ of uncountable coloring number, there is a subgraph $H$ 
  of $G$ of size $\aleph_1$ with uncountable coloring number.
\end{xitemize}
This implication can be proved directly but we can also see this using the terminology 
introduced in  
\cite{I} as follows: The downward Löwenheim-Skolem Theorem \equationof{eq:-0} for
$L_{stat}$ is equivalent to the Diagonal Reflection Principle $\DRP(\IC_{\aleph_0})$ down 
to an internally club set (Corollary 3.6 in \cite{I}). This implies the reflection 
principle $\RP_{\IU_{\aleph_0}}$ down to an internally unbounded set of size $\LT\aleph_2$. 
This reflection 
principle is equivalent to Axiom R of Fleissner (Lemma 2.6 in \cite{RIMS13}). From Axiom 
R, the Fodor-type Reflection Principle (\FRP) follows (Corollary 2.6 in 
\cite{fuchino-juhasz-etal}).  \equationof{eq:-1} is a consequence of (actually equivalent 
to \FRP\ over \ZFC\ (\cite{more})). 

As it is mentioned in the proof of \theoremof{p-0}, the condition ``$\aleph_2$ 
is $\calP$-generic supercompact for $\calP=$ the class of all 
$\LT\aleph_1$-closed \pos.'' is 
equivalent to Game Reflection Principle (\GRP). As the name suggests, \GRP\ 
is actually a reflection principle which claims the reflection of the non-existence of 
winning strategy of certain  
games of length $\omega_1$ down to subgames of size $\LT\aleph_2$. A remarkable feature of 
this principle is that it implies \CH\ (Theorem 8 in König \cite{koenig}).

The notion of Laver-generic large cardinals was introduced in Fuchino, Ottenbreit Maschio 
Rodrigues and Sakai\cite{II} in search for 
reflection principles which generalize \GRP. The following definition of Laver-generic 
large cardinals is a streamlined version adopted in later papers Fuchino and Ottenbreit Maschio 
Rodrigues \cite{sf-omr}, 
Fuchino \cite{future} etc.\ and slightly different from the one given in \cite{II}.

We call a non-empty class $\calP$ of \pos\ {\It iterable} if it satisfies: \Hypertarget{iterable}{}
\ifextended\begin{itemize}\item[\daimaru{1}]\quad\fi
$\ssetof{\bbone}\in\calP$,
\ifextended\quad\item[\daimaru{2}]\ \ \fi
$\calP$ is closed \wrt\ forcing 
equivalence \ifextended\\ \fi
(i.e.\ if $\poP\in\calP$ and $\poP\sim\poP'$ then $\poP'\in\calP$), 
\ifextended\item[\daimaru{3}]\quad
\fi
closed \wrt\ restriction \ifextended\\ \fi
(i.e.\ if $\poP\in\calP$ then $\poP\restr\condp\in\calP$ for any $\condp\in\poP$), and, 
\ifextended\item[\daimaru{4}]\quad
\fi
for 
any $\poP\in\calP$ and $\poP$-name $\utpoQ$, $\forces{\poP}{\utpoQ\in\calP}$ implies
$\poP\ast\utpoQ\in\calP$.\ifextended\end{itemize}\fi

For an iterable class $\calP$ of \pos, a cardinal $\kappa$ is said to 
be {\It $\calP$-Laver-generically supercompact}\/ if, for 
any $\lambda>\kappa$ and $\poP\in\calP$, there is a $\poP$-name $\utpoQ$ with 
$\forces{\poP}{\utpoQ\in\calP}$, \st\  
for $(\uniV,\poP\ast\utpoQ)$-generic $\genH$, there  
are $j, M\subseteq\uniV[\genH]$ 
\st\ \smash{$\Elembed{j}{\uniV}{M}{\kappa}$}, $j(\kappa)>\lambda$, 
$\poP,\genH, j\imageof\lambda\in M$.

$\kappa$ is {\It tightly $\calP$-Laver-generically supercompact} if it is 
$\calP$-Laver-generically supercompact and $\utpoQ$, $j$ and $M$ for each $\poP\in\calP$ in 
the definition of $\calP$-Laver-generic supercompactness additionally satisfy 
that $\poP\ast\utpoQ$ is forcing 
equivalent to a \po\ of size $\LE j(\kappa)$.\footnote{\Label{fn-0} In the 
  following, we shall 
  denote this condition simply by 
  ``$\cardof{\poP\ast\utpoQ}\leq\lambda$''. 
  More generally, we shall simply write ``$\cardof{\poP}\leq\lambda$'' for a \po\ $\poP$ to 
  say that ``the \po\ $\poP$ is forcing equivalent to a \po\ of size $\leq\lambda$".} 

\ifextended{\extendedcolor 
A cardinal $\kappa$ is {\It(tightly) $\calP$-Laver-generically superhuge}, if $\kappa$ 
satisfies the conditions of (tightly) $\calP$-Laver-generically supercompactness, with the 
condition $j\imageof\lambda\in M$ replaced by $j\imageof j(\kappa)\in M$. 
Clearly a 
(tightly) $\calP$-Laver-generically superhuge cardinal is 
(tightly) $\calP$-Laver-generically supercompact. 

The name ``Laver-generic large cardinal'' is chosen in connection with fact that 
Laver-function plays central role in the construction of standard models with Laver-generic 
large cardinals (see \theoremof{p-3}). 

Laver-genericity corresponding to other notions of large cardinals can be defined 
canonically: A cardinal 
$\kappa$ is {\It$\calP$-Laver-generically ultrahuge}, if it enjoys 
the definition of $\calP$-Laver-generically supercompactness and that the condition
``$j\imageof{\lambda}\in M$'' in the definition of supercompactness is replaced by the 
stronger ``${\uniV_{j(\lambda)}}^{\uniV[\genH]}\in M$''.
{\ifextended\extendedcolor\\That is:\Hypertarget{LGHH}{}
\begin{itemize} 
\item[]\hspace{-1.2em}\extendedcolor $\kappa$ is {\It(tightly) $\calP$-Laver-generically 
  ultrahuge}, if 
  for any $\lambda>\kappa$ and $\poP\in\calP$ there is a $\poP$-name $\utpoQ$ with 
  $\forces{\poP}{\utpoQ\in\calP}$, \st\  
  for $(\uniV,\poP\ast\utpoQ)$-generic $\genH$, \\there  
  are $j, M\subseteq\uniV[\genH]$ 
  \st\ 
  $\Elembed{j}{\uniV}{M}{\kappa}$, $j(\kappa)>\lambda$, 
  $\poP,\genH, {V_{j(\lambda)}}^{\uniV[\genH]}\in M$,\\
  \smash{(and
  $\cardof{\poP\ast\utpoQ}\leq j(\kappa)$).} 
\end{itemize}\fi}
$\kappa$ is {\It$\calP$-Laver-generically hyperhuge} if $\kappa$ 
satisfies the  
definition  
obtained by replacing ``$j\imageof\lambda\in M$'' in 
the definition of $\calP$-Laver-generically supercompactness by
``$j\imageof{j(\lambda)}\in M$''.  
\ifextended{\extendedcolor\\That is:
\begin{itemize} 
\item[]\hspace{-1.2em}\extendedcolor $\kappa$ is {\It(tightly) $\calP$-Laver-generically 
  hyperhuge}, if 
  for any $\lambda>\kappa$ and $\poP\in\calP$ there is a $\poP$-name $\utpoQ$ with 
  $\forces{\poP}{\utpoQ\in\calP}$, \st\  
  for $(\uniV,\poP\ast\utpoQ)$-generic $\genH$, \\there  
  are $j, M\subseteq\uniV[\genH]$ 
  \st\ 
  $\Elembed{j}{\uniV}{M}{\kappa}$, $j(\kappa)>\lambda$, 
  $\poP,\genH, j\imageof{j(\lambda)}\in M$,\\
  \smash{(and
  $\cardof{\poP\ast\utpoQ}\leq j(\kappa)$).} 
\end{itemize}}\memo{\normalsize!!!}\smallskip\else
{\extendedcolor``Tight''} versions of these notions are also defined similarly as above. 
\fi
}\else

Laver-generic large cardinals corresponding to other notions of large large cardinals can 
be defined similarly. In particular we consider in the following (tightly) $\calP$-generic 
superhuge/ultrahuge/hyperhuge cardinals which correspond to superhuge/ultrahuge/hyperhuge 
cardinals respectively. For the precise definition of these cardinals see \cite{II}, \cite{future} 
and/or \cite{laver-gen-maximum} or the extended postprint version of the present article 
mentioned in the  
footnote on \pageof{postprint}. 
\fi

The following implications follow from the definitions:
\begin{equation*}
  \begin{array}{c}
    \kappa\mbox{ is (tightly) }\calP\mbox{-Laver-generically ultrahuge}\\[\jot]
    \Downarrow\\[\jot]
    \kappa\mbox{ is (tightly) }\calP\mbox{-Laver-generically superhuge}\\[\jot]
    \Downarrow\\[\jot]
    \kappa\mbox{ is (tightly) }\calP\mbox{-Laver-generically supercompact}
  \end{array}
\end{equation*}

The relationship between Laver-generic hyperhugeness and Laver-generic ultrahugeness is  
slightly more subtle and at the moment, we need the tightness to obtain the expected implication:

\begin{Lemma}{\rm (Fuchino and Usuba \cite{laver-gen-maximum})}
  \Label{p-0-0} For any class $\calP$ of \pos, 
  if $\kappa$ is tightly $\calP$-Laver-generically hyperhuge then $\kappa$ is tightly
  $\calP$-Laver-generically ultrahuge. \ifextended\else\qed\fi
\end{Lemma}
{\ifextended\extendedcolor
\prf Suppose that $\kappa$ is tightly $\calP$-Laver-generically hyperhuge, and  
$\lambda>\kappa$. \Wolog, we may assume that 
\begin{equation}\Label{x-a-a-0}
  V_\lambda\prec_{\Sigma_n}\uniV\mbox{ for 
  sufficiently large }n. 
\end{equation}
Let $\lambda^*:=(\cardof{V_\lambda}^+)^\uniV$. For $\poP\in\calP$, let $\utpoQ$ be 
a $\poP$-name \st\ $\forces{\poP}{\utpoQ\in\calP}$ and, 
for $(\uniV,\poP\ast\utpoQ)$-generic $\genH$, there are $j$, $M\subseteq\uniV[\genH]$ \st\
\begin{equation}\Label{x-a-0}
  \Elembed{j}{\uniV}{M}{\kappa},\ j(\kappa)\geq\lambda^*, 
\end{equation}
\begin{equation}\Label{x-0}
  j\imageof{j(\lambda^*)},\,\poP,\,\genH\in M,\mbox{ and}
\end{equation}
\begin{equation}\Label{x-1}
  \cardof{\poP\ast\utpoQ}\leq j(\kappa). 
\end{equation}
\begin{Claim}\Label{claim-0}
  For $\alpha\leq j(\lambda)$, ${V_\alpha}^\uniV\in M$.
\end{Claim}
\prfofClaim
By induction on $\alpha\leq j(\lambda)$, we prove
\begin{equation}\Label{x-2}
  {V_\alpha}^\uniV\in M\mbox{ and }{V_\alpha}^\uniV\subseteq {V_\alpha}^M.
\end{equation}

For $\alpha<\omega$, this is clear. Suppose we have
${V_\alpha}^\uniV\in M$ and ${V_\alpha}^\uniV\subseteq {V_\alpha}^M$. 
Then, since $M\models\cardof{{V_\alpha}^M}<j(\lambda^*)$ (by the choice 
of $\lambda^*$ and) by elementarity,
$\calP^\uniV({V_\alpha}^\uniV)\subseteq\calP^M({V_\alpha}^M)\subseteq M$ 
by \equationof{x-0} and Lemma 2.5,\,\assertof{5} in \cite{II}. Again by Lemma 
2.5,\,\assertof{5} in \cite{II}, it follows that
${V_{\alpha+1}}^\uniV=\calP^\uniV({V_\alpha}^\uniV)\in M$, and 
${V_{\alpha+1}}^\uniV=\calP^\uniV({V_\alpha}^\uniV)
\subseteq\calP^M({V_\alpha}^M)={V_{\alpha+1}}^M$ 

If $\gamma\leq j(\lambda)$ is a limit, and $V_\alpha\in M$, 
$M\models{V_\alpha}^\uniV\subseteq V_\alpha$ for all $\alpha<\gamma$, then
$\seqof{{V_\alpha}^\uniV}{\alpha<\gamma}\subseteq M$. Hence by Lemma 2.5,\,\assertof{5} in 
\cite{II}, it follows that $\seqof{{V_\alpha}^\uniV}{\alpha<\gamma}\in M$. Thus 
${V_\gamma}^\uniV=\bigcup_{\alpha<\gamma}{V_\alpha}^\uniV\in M$ and 
${V_\gamma}^\uniV=\bigcup_{\alpha<\gamma}{V_\alpha}^\uniV\subseteq
\bigcup_{\alpha<\gamma}{V_\alpha}^M={V_\gamma}^M$. 
\qedofClaim\qedskip

Now, it follows that 
\begin{equation*}
  M\ubecause{\ni}{}{by \Claimof{claim-0} and \equationof{x-0}}
  V_{j(\lambda)}[\genH]\obecause{=}{}{by \equationof{x-a-a-0} and \Lemmaof{p-Lg-RcA-0-0}} 
  {V_{j(\lambda)}}^{\uniV[\genH]}.
\end{equation*}

This shows that $j$ and $M$ taken here are as in the definition 
of $\calP$-Laver-generically ultrahugeness. 
\qedofLemma
\qedskip\fi}

The consistency strength of Laver-generic hyperhugeness can be separated from that of other 
notions of 
Laver-generic large cardinals. This is because we know that the existence of $\calP$-Laver-generic 
hyperhuge cardinal for any $\calP$ (if it ever exists) is equiconsistent with that 
of a hyperhuge cardinal (see \theoremof{p-bedrock-2}).

\memox{\normalsize !!!}
At first glance,
it is not immediately clear if the notion of Laver-generic large 
cardinal is definable in the language $\Lin$ of \ZFC. In \cite{definability} 
an abstract generic version of extender is introduced to show the definability of 
Laver-generic large cardinals.

Laver-generic supercompactness implies double plus versions of forcing axioms.
For a class $\calP$ of \pos\ and cardinals $\kappa$, $\mu$, we denote with 
${\sf MA}^{+\mu}(\calP,\LT\kappa)$ and ${\sf MA}^{++\LE\mu}(\calP,\LT\kappa)$ the 
following versions of Martin's Axiom:

\begin{itemize}
\item[\darkred$\darkred{\sf MA}^{+\mu}(\calP,\LT\kappa)$: ]\qquad\qquad\qquad
  For any $\poP\in\calP$, 
  any family  
  $\calD$ of dense subsets of\/ $\poP$ with $\cardof{\calD}<\kappa$ and any family 
  $\calS$ of\/ $\poP$-names \st\ $\cardof{\calS}\leq\mu$ and
  $\forces{\poP}{\utilde{S}\xmbox{ is a stationary subset of }\omega_1}$ for all
  $\utilde{S}\in\calS$, there is a $\calD$-generic filter $\genG$ over $\poP$ \st\
  $\utilde{S}[\genG]$ is a stationary subset of $\omega_1$ for all
  $\utilde{S}\in\calS$. 
\end{itemize}\begin{itemize}
\item[\darkred$\darkred{\sf MA}^{++\LE\mu}(\calP,\LT\kappa)$: ]\qquad\qquad\qquad\quad
  For any $\poP\in\calP$, any family 
  $\calD$ of dense subsets of\/ $\poP$ with $\cardof{\calD}<\kappa$ and any family 
  $\calS$ of\/ $\poP$-names \st\ $\cardof{\calS}\leq\mu$ and
  \smash{$\forces{\poP}{\utilde{S}\xmbox{ is a stationary subset of }
      \Pkl{\eta_{\scriptstyle\utilde{S}}}{\theta_{\utilde{S}}}}$} 
  for some $\omega<\eta_{\utilde{S}}\leq\theta_{\utilde{S}}\leq\mu$ 
  with $\eta_{\utilde{S}}$ regular, for all 
  $\utilde{S}\in\calS$, there is a $\calD$-generic filter $\genG$ over $\poP$ \st\
  $\utilde{S}[\genG]$ is stationary in
  $\Pkl{\eta_{\scriptstyle\utilde{S}}}{\theta_{\utilde{S}}}$ for all 
  $\utilde{S}\in\calS$.  
\end{itemize}

Clearly $\MA^{++\LT\omega_2}(\calP,\LT\kappa)$ is 
equivalent to $\MA^{+\omega_1}(\calP,\LT\kappa)$.\\
$\MM^{++}$ is $\MA^{+\omega_1}(\mbox{stationary preserving \pos},\,\LT\aleph_2)$. 

\begin{theorem}{\rm (Theorem 5.7 in Fuchino, Ottenbreit Maschio 
Rodrigues and Sakai \cite{II}, see also Fuchino \cite{future})}
  \Label{p-1} For an iterable class $\calP$ of \pos\ \st\ 
  \begin{xitemize}
  \xitem[x-laver-6] 
    the elements of $\calP$ preserve stationarity of subsets of $\Pkl{\mu}{\theta}$ for all 
    $\mu\leq\theta<\kappa$,
  \end{xitemize}
  if $\kappa>\aleph_1$  
  is $\calP$-Laver-generically supercompact then $\MA^{++\leq\mu}(\calP,\LT\kappa)$ holds 
  for all 
  $\mu<\kappa$. \qed
\end{theorem}

In contrast to usual generic large cardinals, a Laver-generic large cardinal if it 
exists, is unique and it is $\kappa_\refl$ in many cases. 

\begin{Lemma}\Label{p-Lg-RA-1-2-0}{\rm(\cite{future}, \cite{laver-gen-maximum}, see also 
    Proposition 4, in \cite{nagoya})}  
  \wassertof{1} If 
  $\kappa$ is $\calP$-generically measurable for an $\omega_1$ preserving iterable 
  $\calP$, then $\omega_1<\kappa$.\smallskip

  \wassert{2} If $\kappa$ is $\calP$-Laver-generically supercompact for an 
  $\omega_1$-preserving iterable $\calP$ with 
  $\Col(\omega_1,\ssetof{\omega_2})\in\calP$ then $\kappa=\omega_2$.\smallskip

  \wassert{3} If $\kappa$ is $\calP$-Laver-generically supercompact for an iterable $\calP$ 
  which contains a \po\ adding a new real, then $\kappa\leq2^{\aleph_0}$.\smallskip

  \wassert{4} If $\kappa$ is $\calP$-generically supercompact for an iterable $\calP$ \st\ 
  all \pos\ in 
  $\calP$ do not add any reals then $2^{\aleph_0}<\kappa$.\smallskip

  \wassert{5} If $\kappa$ is $\calP$-Laver-generically supercompact for an iterable $\calP$ 
  which contains a \po\ which collapses $\aleph_1$ then $\kappa=\aleph_1$.  
\end{Lemma}
\prf We only prove \assertof{5} since it is not explicitly given in 
\cite{laver-gen-maximum}. Suppose that $\kappa$ is $\calP$-Laver-generically supercompact 
and $\poP\in\calP$ is \st\ $\forces{\poP}{{\aleph_1}^\uniV\mbox{ is countable}}$. If
$\kappa\not=\aleph_1$, then we have $\aleph_1<\kappa$. Let $\utpoQ$ be a $\poP$-name of a 
\po\ \st, for $(\uniV,\poP\ast\utpoQ)$-generic $\genH$, there are $j$,
$M\subseteq\uniV[\genH]$ \st\ $\Elembed{j}{\uniV}{M}{\kappa}$ and $\poP$, $\genH\in M$. By
$\genH\cap\poP\in M$ and since ${\aleph_1}^\uniV<\crit(j)$, we have
$M\modelof{{\aleph_1}^\uniV=j({\aleph_1}^\uniV)\mbox{ is countable}}$. This is a 
contradiction to the elementarity of $j$. \qedofLemma\qedskip

A cardinal $\kappa$ is called {\It greatly weakly Mahlo} if $\kappa$ is weakly inaccessible 
and 
there exists a non-trivial $\LT\kappa$-complete normal filter $\calF$ over $\kappa$ \st\
$\setof{\mu<\kappa}{\mu\xmbox{ is a regular cardinal}}\in\calF$, and $\calF$ is closed 
\wrt\ the Mahlo operation $\mahlo$ \footnote{Closedness here means that for 
  any $S\in\calF$, we have 
  $\mahlo(S)\in\calF$.} where 
\begin{itemize}
\item[]
  $S\ \mapsto\ {\darkred\mahlo(S)}:=\setof{\alpha\in S}{{}
  \begin{array}[t]{@{}l}
    \alpha\mbox{ has uncountable cofinality and}\\
    S\cap\alpha\mbox{ is stationary in }\alpha} \mbox{\qquad\qquad(\cite{ccc}).}\end{array}$
\end{itemize}
Note that the Mahlo operation given above is slightly different from the one in 
\cite{baumgartner-taylor-wagon}. 

If $\kappa$ is greatly weakly Mahlo then it is hyper-weakly Mahlo (Proposition 3.4 in 
\cite{ccc}). 

The tightness of the Laver-genericity can be still strengthened as follows: a cardinal 
$\kappa$ is {\It tightly$^+$ $\calP$-Laver-generically $x$-large}, for a 
notion ``$x$-large'' of large cardinal (e.g.\ ``supercompact'', ``superhuge'' 
etc.) if it satisfies the definition of tightly $\calP$-Laver-generically $x$-large 
cardinal with the condition ``$\cardof{\poP\ast\utpoQ}\leq j(\kappa)$'' in the definition 
being replaced by the condition ``there is a complete Boolean algebra $\BaB$ of size
$j(\kappa)$ \st\ $\BaB^+$ is forcing equivalent to $\poP\ast\utpoQ$''. 

\begin{theorem}\Label{p-1-0} \wassertof{1} {\rm (Theorem 3.5 in \cite{ccc})} If $\kappa$ is a
  $\ssetof{\poP}$-generically measurable 
  for a \po\ $\poP$ with the $\mu$-cc for some $\mu<\kappa$, then $\kappa$ is greatly 
  weakly Mahlo. \smallskip

  \wassert{2} {\rm (Theorem 5.8 in Fuchino, Ottenbreit Maschio 
Rodrigues and Sakai \cite{II})} If $\kappa$ is tightly $\calP$-Laver-generically superhuge for a class $\calP$ of ccc \pos\ \st\ at least one element 
  of $\calP$ adds a 
  real, then $\kappa=\continuum$.\smallskip

  \wassert{3} For an iterable class $\calP$ of \pos, if 
  $\kappa$ is tightly$^+$ $\calP$-Laver-generically hyperhuge, then $2^{\aleph_0}\leq\kappa$. 
  \qed 
\end{theorem}

We give a sketch of the proof of \theoremof{p-1-0},\,\assertof{3} after 
\Corof{p-bedrock-3}
\memox{{\normalsize!!!} recurrence-axioms.tex Lemmaof p-Lg-RcA-4-1}. 

In Fuchino and Usuba \cite{laver-gen-maximum}, it is proved as a Corollary of \theoremof{p-12} 
and \theoremof{p-bedrock-2} that for $\kappa=\aleph_1$, the `$+$' in ``tightly$^+$'' in 
\theoremof{p-1-0},\,\assertof{3} above can be dropped. 

\begin{theorem}\Label{p-2}{\rm(Fuchino, Ottenbreit Maschio 
Rodrigues and Sakai \cite{II}, Fuchino and Usuba \cite{laver-gen-maximum} for \assertof{$\Delta$} 
    and \assertof{$\Delta'$})}
  \wassertof{$A$} If $\calP$ is the class of all $\LT\aleph_1$-closed \pos,  
  and $\kappa$ is 
  $\calP$-Laver-generically supercompact, then $\kappa=\aleph_2$ and \CH\ holds.\smallskip

  \wassert{$B$} If $\calP$ is either the class of all proper \pos\ or the class of all 
  semi-proper \pos, and $\kappa$ is $\calP$-Laver-generically supercompact, then
  $\kappa=\continuum=\aleph_2$.\smallskip 

  \wassert{$\Gamma$} If $\calP$ is the class of all ccc \pos, and $\kappa$ is 
  $\calP$-Laver-generically supercompact, then $\kappa$ is very large 
  and $\kappa\leq\continuum$. \smallskip 

  \wassert{$\Gamma'$} If $\calP$ is the class of all ccc \pos, and $\kappa$ is 
  tightly $\calP$-Laver-generically superhuge, then $\kappa$ is very large and
  $\kappa=\continuum$.\smallskip

  \wassert{$\Delta$} If $\calP$ is the class of all \pos, and $\kappa$ is $\calP$-Laver-generically supercompact, then $\kappa=\aleph_1$. 

  \wassert{$\Delta'$} If $\calP$ is the class of all \pos, and $\kappa$ is tightly$^+$
  $\calP$-Laver-generically supercompact, then $\kappa=\aleph_1$. 
\end{theorem}
\prf \assertof{$A$}: By \Lemmaof{p-Lg-RA-1-2-0}, \assertof{2},\,\assertof{4}.\quad
\assertof{$B$}: By \Lemmaof{p-Lg-RA-1-2-0}, \assertof{2},\,\assertof{3} and \theoremof{p-1}.
\quad \assertof{$\Gamma$}: By \Lemmaof{p-Lg-RA-1-2-0}, \assertof{3}, and 
\theoremof{p-1-0},\,\assertof{1}. \quad \assertof{$\Gamma'$}: By \assertof{$\Gamma$} and 
\theoremof{p-1-0},\,\assertof{2}.\\ \assertof{$\Delta$}: By \Lemmaof{p-Lg-RA-1-2-0}, 
\assertof{5}. \quad \assertof{$\Delta'$}: By \assertof{$\Delta$} and 
\theoremof{p-1-0},\,\assertof{3}. 
\qedoftheorem\qedskip

The consistency of the existence of a $\calP$-Laver-generic large cardinal can be proved 
under the existence of corresponding genuine large cardinal.

\begin{theorem}{\rm (Theorem 5.2, \cite{II})} 
  \Label{p-3} \wassertof{$A$} Suppose that $\kappa$ is supercompact
  and $\poP=\Col(\aleph_1,\kappa)$, then, 
  in $\uniV[\genG]$, for any $(\uniV,\poP)$-generic $\genG$, ${\aleph_2}^{\uniV[\genG]}$
  ($=\kappa$) is tightly 
  $\calP$-closed-Laver-generically supercompact for the class $\calP$ of 
  all $\sigma$-closed \pos\ (and $\CH$ holds).\smallskip

  \wassert{$B$}
  Suppose that $\kappa$ is superhuge with a
  Laver function $\mapping{f}{\kappa}{V_\kappa}$ for superhugeness.  If $\poP$ is the 
  CS-iteration for forcing \PFA\ along with $f$, then, 
  in $\uniV[\genG]$ for any $(\uniV,\poP)$-generic $\genG$, ${\aleph_2}^{\uniV[\genG]}$
  ($=\kappa$)  
  is tightly$^+$ $\poP$-Laver-generically superhuge for the class $\calP$ of all proper 
  \pos\ (and $\continuum=\aleph_2$ holds).\smallskip

  \wassert{{$B'$}}
  Suppose that $\kappa$ is superhuge with a
  Laver function $\mapping{f}{\kappa}{V_\kappa}$ for superhugeness. If
  $\poP$ is the RCS-iteration for forcing \MM\ along with $f$, then,  
  in $\uniV[\genG]$ for any $(\uniV,\poP)$-generic $\genG$, ${\aleph_2}^{\uniV[\genG]}$
  ($=\kappa$)  
  is tightly$^+$ $\calP$-Laver-generically superhuge for the class $\calP$ of all 
  semi-proper \pos\ 
  (and 
  $\continuum=\aleph_2$ holds).\smallskip

  \wassert{$\Gamma$}
  Suppose that $\kappa$ is supercompact with a 
  Laver function $\mapping{f}{\kappa}{V_\kappa}$ for supercompactness. If $\poP$ is a 
  FS-iteration for forcing \MA\ along with $f$,   
  then, in $\uniV[\genG]$ for any $(\uniV,\poP)$-generic $\genG$, $\continuum$ ($=\kappa$) is 
  tightly$^+$ $\calP$-Laver-generically supercompact for the class $\calP$ of all ccc \pos 
  (and $\kappa=\continuum$ while $\kappa$ still is 
  very large).\smallskip

  \wassert{$\Delta$} 
  Suppose that $\kappa$ is supercompact with a 
  Laver function $\mapping{f}{\kappa}{V_\kappa}$ for supercompactness. If $\poP$ is a 
  FS-iteration for forcing $f$ where $f$ is used to book-keep through all \pos\ 
  in $V_\kappa$, then in $\uniV[\genG]$ for 
  any $(\uniV,\poP)$-generic $\genG$, $\continuum$ ($=\kappa$) is  
  tightly$^+$ $\calP$-Laver-generically supercompact for the class $\calP$ of all \pos\  
  (and $\CH$ holds).\qed
\end{theorem}

\theoremof{p-3} above also holds for all other notions of large cardinal and corresponding 
Laver-generic version of generic large cardinal except \assertof{$B$} and \assertof{$B'$} 
in which the supercompactness does not seem to be strong enough to show that the resulting 
generic extension in the proof satisfies  
the expected Laver-genericity. 

In a sense, the cases treated in \theoremof{p-3} are (almost) exhaustive. This can be seen 
in the following:

\begin{theorem}\Label{p-3-0}
  Suppose that
  $\calP$ is an iterable class of \pos\ \st\ all\/ $\poP\in\calP$ 
  are $\omega_1$-preserving and $\calP$ contains a \po\/ $\poP^*$ whose generic filter 
  destroys a stationary subset of $\omega_1$.\footnotemark\ Then there is no 
  $\calP$-Laver-generically supercompact cardinal.
\end{theorem}
\prf\footnotetext{``$\poP^*$ destroys a stationary subset of $\omega_1$'' means here that a
  $\poP^*$-generic set codes a club subset of $\omega_1\setminus S$ in some absolute way. 
  \par 
  Note that, for stationary and co-stationary subset $S$ of $\omega_1$, various 
  \pos\ are known which preserve $\omega_1$ while shooting a club in
  $\omega_1\setminus S$\ifextended\ (e.g.\ see \cite{abraham-shelah}).\else.\fi}%
Suppose, toward a contradiction, that $\calP$ is as above and there 
is $\calP$-Laver-generically supercompact cardinal $\kappa$. 

Let $S\subseteq\omega_1$ be 
stationary \st\ there is a \po\ $\poP^*\in\calP$ shooting a club in
$\omega_1\setminus S$. Let $\lambda>\cardof{\poP^*}$ be large enough. 
By assumption, there is 
a $\poP^*$-name $\utpoQ$ of a \po\ \st\ $\forces{\poP^*}{\utpoQ\in\calP}$ and, for
$(\uniV,\poP^*\ast\utpoQ)$-generic $\genH$, there are $j$, $M\subseteq\uniV[\genH]$ \st\ 
\begin{equation}
  \Label{x-Lg-RA-0} 
  \Elembed{j}{\uniV}{M}{\kappa},
\end{equation}
\begin{equation}
  \Label{x-Lg-RA-0-0}
  j(\kappa)>\lambda,\quad\mbox{and}
\end{equation}
\begin{equation}
  \Label{x-Lg-RA-1} j\imageof{\lambda}, \poP,\ \genH\in M.
\end{equation}

By the choice of $\poP^*$ ($\circleq\poP^*\ast\utpoQ$), 
there is a nice $\calP^*$-name $\utilde{C}\in\uniV$ of a club set 
$\subseteq{\omega_1}^\uniV\setminus S$. By \xitemof{x-Lg-RA-0-0}, \xitemof{x-Lg-RA-1} and 
by the choice of $\lambda$, we have $\utilde{C}\in M$. 

Thus 
$M\modelof{S\mbox{ is a non-stationary subset of }\omega_1}$ by \xitemof{x-Lg-RA-1}. Since
$\crit(j)=\kappa>\omega_1$ by \Lemmaof{p-Lg-RA-1-2-0},\,\assertof{1}\imemox{!!!}, 
we have $S=j(S)$. By $\uniV\modelof{S\xmbox{ is stationary subset of }\omega_1}$, this is a 
contradiction to the elementarity \equationof{x-Lg-RA-0} of $j$.  
\qedoftheorem
\qedskip

$\calP$-Laver genericity for stationary preserving $\calP$, in particular those $\calP$ 
containing all $\sigma$-closed \pos\ can be regarded as a strong reflection principle.  

\begin{theorem}
  \Label{p-4} Suppose that $\calP$ is an iterable class of \pos\ which 
  are $\omega_1$-preserving and include 
  all $\sigma$-closed \pos. If $\kappa$ is $\calP$-Laver-generically supercompact then 
  {
    \begin{xitemize}
      \xitemas[eq:-0] 
	    for any structure $\gmA$ of countable signature, there 
      is $\gmB\prec_{L_{stat}}\!\gmA$
      of cardinality $\LT\aleph_2$.
  \end{xitemize}}
\end{theorem}
\prf 
By \theoremof{p-1}, $\MA^{++}(\sigma\mbox{-closed})$ holds. Cox \cite{cox:} proved that
$\MA^{++}(\sigma\mbox{-closed})$ implies $\DRP(\LT\aleph_2,\IC_{\aleph_0})$ (in the 
notation of \cite{I}). By Lemma 3.5 in \cite{I}, this principle is equivalent to 
\equationof{eq:-0}. \qedoftheorem\qedskip

The existence of $\calP$-Laver generic large cardinal for the class of all ccc-\pos\ also 
implies a reflection statement similar to \equationof{eq:-0} (see Theorem 
5.9,\,\assertof{3} in \cite{II}). It is not known if the existence of $\calP$-Laver generic 
large cardinal for the class of all ccc-\pos\ imply \FRP\ but it easy to see that we can 
force \FRP\ and the existence of the Laver generic large cardinal starting from a model 
with a supercompact with the relevant large 
cardinal above it (Fuchino, Ottenbreit Maschio Rodrigues, and 
Sakai \cite{II}). \qedskip

In summary, Laver-generic large cardinal axiom (i.e.\ the axiom asserting the 
existence of a $\calP$-Laver-generic large cardinal for some iterable class of \pos\ $\calP$) is 
an assertion about the existence of a 
$\calP$-generic large cardinal (without "Laver-" !) with the feature of resurrection 
property similar to the one in 
Resurrection Axioms studied by Hamkins and Johnstone \cite{hamkins-johnstone}, 
\cite{hamkins-johnstone2}, and Tsaprounis \cite{tsaprounis1}: for any generic extension by 
some $\poP\in\calP$ we find a further generic extension by a \po\ of the form
\smash{$\poP\ast\utpoQ\in\calP$} in which a previously chosen instance of generic elementary 
embedding corresponding to the generic largeness of the cardinal resurrects. 

In spite of this similarity, it is only recent that we found there is a connection 
between these two  
similar types of axioms. First remember that, for each reasonable class $\calP$ consisting 
of proper or stationary preserving 
\pos\ $\calP$, $\calP$-Laver generic large cardinal axiom implies the $\calP$-Laver generic 
large cardinal is $\kappa_\refl$ (\theoremof{p-3}). 
In Fuchino \cite{future} it is then proved that the tightly $\calP$-Laver-generic superhuge  
axiom implies the boldface version Resurrection Axiom (in the sense of Hamkins and 
Johnstone) for parameters $\calP$ and  
$\calH(\kappa_\refl)$ (Fuchino \cite{future}).
The tightly $\calP$-Laver generic ultrahuge axiom implies the tight 
version of the Unbounded Resurrection Axiom of Tsaprounis for $\calP$ (Fuchino \cite{future}). 

$\calP$-Laver generic large cardinal axioms for reasonable classes $\calP$ consisting 
of proper or stationary preserving \pos\ imply double plus versions of Forcing Axiom 
(\theoremof{p-1}). By \theoremof{p-0} (and its proof), \theoremof{p-4} and the remarks 
after that most of the known structural 
reflection principles down to $\LT\kappa_\refl$ are covered by these axioms. 

$\calP$-Laver generic large cardinal axioms also support the intuition shared by 
many set-theorists that the continuum is either $\aleph_1$ or $\aleph_2$ or very large 
(\theoremof{p-2}, \theoremof{p-3} and \theoremof{p-3-0}). Perhaps we can also formulate 
this circumstance as: if we are to  
support one of the reasonable instances of Laver-generic large cardinal axiom then the 
continuum must be either $\aleph_1$ or $\aleph_2$ or extremely large ($\aleph_{\continuum}$ 
in particular, and much more). 

Laver generic large cardinal axiom can be considered as a strong reflection principle  
with the feature of resurrection (à la Hamkins and Johnstone). Certain amount of 
absoluteness is also inherent in Laver generic large cardinal axiom in that it implies 
strong versions of forcing axioms (\theoremof{p-1}). 
In \sectionof{sec:3} we shall see 
that more generic absoluteness in terms of recurrence axioms (see the next section for 
these axioms) can be integrated into the 
Laver genericity by moving to the notions of large cardinals with the additional 
properties (in plural since they can be formulated only in an axiom scheme --- or in some cases 
even not first-order formalizable) which we call "super $C^{(\infty)}$" 
(\theoremof{p-Lg-RcA-5}) and that this is still realizable below the consistency strength 
of 2-huge (\theoremof{p-Lg-RcA-2} and \theoremof{p-Lg-RcA-4}).

We can also feel more confident (or insecure in case you do not believe in large large 
cardinals\footnote{Note (not only) for the editor: ``large large cardinals'' is not a 
  typo.}) with these stronger versions of Laver-generic large cardinal axiom since their  
consistency strength can be decided (see Corollaries \ref{p-bedrock-5}, \ref{p-bedrock-6}, 
\ref{p-bedrock-7}, \ref{p-bedrock-8}). 

\section{Recurrence, Maximality, and the solution(s) of the Continuum Problem}
\Label{sec:2}
The Recurrence Axioms we are going to introduce below can be considered as reflection 
statements in 
set-generic multiverse down to set-generic geology. Their strengthenings can be also considered as 
a sort of absoluteness statements (see \Propof{p-intro-0}). 

Adopting the terminology of set-theoretic geology (see e.g.\ Fuchs, Hamkins, and Reitz \cite{fhr}), 
An inner model $\uniW$ of a universe $\uniU$ (in most of the 
cases $\uniU$ is the 
real universe $\uniV$ but sometimes it is some other universe obtained from $\uniV$) is 
called  
a {\It ground of\/ $\uniU$}, if there is a \po\ $\poP\in\uniW$ and $(\uniW,\poP)$-generic
$\genG\in\uniU$ \st\ $\uniU=\uniW[\genG]$. 

For a class $\calP$ of \pos\ and a set $A$ (of parameters), the {\It Recurrence Axiom for 
$\calP$ and $A$} ({\darkred$(\calP,A)$-\RcA}, for short\footnote{The notation ``\RcA'' is 
  chosen to 
  avoid the collision with ``{\sf RCA}'' which is used in Reverse Mathematics to 
  denote ``recursive comprehension axiom''.}) is the following assertion formulated as an  
axiom scheme in the language $\Lin$ of set theory:

\begin{itemize}
\item[\darkred$(\calP,A)$-\RcA: ]\qquad\qquad\quad For any $\Lin$-formula
  $\varphi=\varphi(\overline{x})$ and 
  $\overline{a}\in A$, if $\forces{\poP}{\varphi(\overline{a}\checked)}$ for 
  a $\poP\in\calP$,  
  then there is a ground $\uniW$ of the universe $V$ \st\ $\overline{a}\in\uniW$ and
  $\uniW\models\varphi(\overline{a})$\footnotemark. 
\end{itemize}
\footnotetext{For sets $a_0$\ctentenc $a_{n-1}$ and $\ol{a}:=\pairof{a_0\ctentenc a_{n-1}}$, we 
  simply write ``$\ol{a}\in A$'' for ``$a_0\ctentenc a_{n-1}\in A$''. For $\poP$-check 
  names (standard $\poP$-names of the elements $a_0\ctentenc a_{n-1}$ in $\uniV$)
  $\check{a}_0\ctentenc\check{a}_{n-1}$ we write
  $\ol{a}\checked:=\pairof{\check{a}_0\ctentenc\check{a}_{n-1}}$.}  

The term ``Recurrence Axiom'' is chosen in allusion to, but not necessarily in (full) agreement with, 
\memox{Nietzsche のスペルを間違えていた}
Nietzsche's \glqq ewige Wiederkehr des Gleichen\grqq\ (eternal recurrence of the same), or 
perhaps rather to Poincaré recurrence theorem: if we understand the relation ``$N$ is 
(set) generic extension of $M$'' as the timeline in the set generic multiverse, we can 
interpret $(\calP,A)$-\RcA\ as saying that 
\begin{xitemize}
\item[] \it if something (formulatable with parameters 
    in $A$) happens in one of the near future universe (in terms of $\calP$) then is is already 
    happened in a not very far past universe (not very far, in the sense that the ``present'' is 
    attainable from there by a set forcing). 
\end{xitemize}
The following is a natural strengthening of the Recurrence Axiom: 
\begin{itemize}
\item[\darkred$(\calP,A)$-\RcAp: ]\qquad\qquad\quad For any $\Lin$-formula
  $\varphi=\varphi(\overline{x})$ and 
  $\overline{a}\in A$, if $\forces{\poP}{\varphi(\overline{a}\checked)}$ for a $\poP\in\calP$, 
  then there is a $\calP$-ground $\uniW$ of the universe $V$ \st\ $\overline{a}\in\uniW$ and
  $\uniW\models\varphi(\overline{a})$. 
\end{itemize}
Here an inner model $\uniW$ of $\uniV$ is called a {\It $\calP$-ground}\/ if there is a 
\po\ $\poP\in\uniW$ with $\uniW\models\poP\in\calP$ and $(\uniW,\poP)$-generic 
$\genG\in\uniV$ \st\ $\uniV=\uniW[\genG]$.

When we regard reflection principles including Laver-generic large cardinal axioms as 
natural axioms, the idea behind this should be that the continuum must be rich enough so 
that many aspects of the situations with sets of high cardinality should be reflected down 
to cardinality $\LT\kappa_\refl$.

Similarly we can consider Recurrence Axioms and their strengthenings as natural 
requirements: 
since our set-theoretic universe $\uniV$ must be highly saturated, as many aspects as 
possible of the situations in set 
generic multiverse should be reflected down to the set generic geology at which we can look 
down form our universe $\uniV$.

We shall use the following version of Laver-Woodin Theorem often without mention. 
It implies in particular that 
$(\calP,A)$-\RcA\ and $(\calP,A)$-\RcAp\ are actually formalizable as 
axiom schemes in $\Lin$.  

\begin{theorem}\Label{p-intro-1}{\rm(Reitz \cite{reitz}, Fuchs-Hamkins-Reitz \cite{fhr})}
  There is an $\Lin$-formula $\Phi(x,r)$ (without any parameters) \st\ the following is provable in \ZFC:\smallskip
  \begin{xitemize}\memoL{\normalsize!!!}
  \xitem[x-intro-10] for all $r$, $\Phi(\cdot,r):=\setof{x}{\Phi(s,r)}$ is a ground 
    in $\uniV$,\smallskip
  \xitem[x-intro-11] for any ground $\uniW$ (of $\uniV$), there is $r$ 
    \st\ $\uniW=\Phi(\cdot,r)$, and\smallskip
  \xitem[x-intro-12]
    if\/ $\uniW$ is a ground of\/ $\uniV$ and $\uniV=\uniW[\genG]$
    where $\genG$ is a $(\uniW,\poP)$-generic 
    for $\poP\in\uniW$, then $r$ \st\ $\uniW=\Phi(\cdot,r)$
    can be chosen as an element of 
    $\psof{(\cardof{\poP}^+)^\uniV}$.
  \end{xitemize}
\end{theorem}

We put together here some other basic facts which will be used in the following.  The next 
lemma was actually already used in the proof of \Lemmaof{p-0-0}. 

\begin{Lemma}{\rm(see Fuchino and Usuba \cite{laver-gen-maximum} for a 
    proof)}\Label{p-Lg-RcA-0-0} If $\alpha$  
  is a limit ordinal and $V_\alpha$ satisfies a  
  large enough fragment of \ZFC, then for any $\poP\in V_\alpha$ and $(\uniV,\poP)$-generic 
$\genG$ we have $V_\alpha[\genG]={V_\alpha}^{\uniV[\genG]}$. \ifextended\else\qed\fi
\end{Lemma}
{\ifextended\extendedcolor
\prf ``$\subseteq$'': This inclusion holds without the condition  
on the fragment of \ZFC. Also the condition ``$\poP\in V_\alpha$'' is irrelevant for this 
inclusion. 

We show by induction on $\alpha\in\On$ that
$V_\alpha[\genG]\subseteq {V_\alpha}^{\uniV[\genG]}$ holds for all $\alpha\in\On$.

The induction steps for $\alpha=0$ and limit ordinals $\alpha$ are trivial. So we assume that
$V_\alpha[\genG]\subseteq {V_\alpha}^{\uniV[\genG]}$ holds and show that the same inclusion 
holds for $\alpha+1$. Suppose $a\in V_{\alpha+1}[\genG]$. Then $a=\uta^\genG$ for 
a $\poP$-name $\uta\in V_{\alpha+1}$. Since $\uta\subseteq V_\alpha$, each
$\pairof{\utb,\condp}\in\uta$ is an element of $V_\alpha$. By induction hypothesis, it 
follows that $\utb^\genG\in {V_\alpha}^{\uniV[\genG]}$. It follows that 
$\uta^\genG\subseteq{V_\alpha}^{\uniV[\genG]}$. Thus 
$a=\uta^\genG\in{V_{\alpha+1}}^{\uniV[\genG]}$. \smallskip 

``$\supseteq$'':  
Suppose that $a\in {V_\alpha}^{\uniV[\genG]}$. Note that we can choose the ``large enough 
fragment of \ZFC'' which should be satisfied in $V_\alpha$ \st\ 
{\ifextended\extendedcolor $(*)$\fi}
${V_\alpha}^{\uniV[\genG]}$ still satisfies a large enough fragment of \ZFC, although the 
fragment may be different from the one $V_\alpha$ satisfies. In particular we find a 
cardinal $\lambda>\cardof{\poP}$ in ${V_\alpha}^{\uniV[\genG]}$ (and hence also in
$\uniV[\genG]$) \st\
$a\in\calH(\lambda)^{V_\alpha^{\uniV[\genG]}}
\subseteq\calH(\lambda)^{\uniV[\genG]}\subseteq {V_\alpha}^{\uniV[\genG]}$. 
Note that\\
    $\calH(\lambda)^{{V_\alpha}^{\uniV[\genG]}}
    =\setof{a}{\cardof{\trcl(a)}<\lambda}^{{V_\alpha}^{\uniV[\genG]}}
    \subseteq\setof{a}{\cardof{\trcl(a)}<\lambda}^{\uniV[\genG]}
    =\calH(\lambda)^{\uniV[\genG]}.$ 

Let
$a^*\in\calH(\lambda)^{\uniV[\genG]}$ be a transitive set \st\ $a\in a^*$. Then $a^*$ 
can be coded by a subset of $\lambda$. We can find the subset of $\lambda$ 
in $\uniV[\genG]$ and this subset has a nice $\poP$-name which is an element of
${V_\alpha}^\uniV$ since $\poP\in V_\alpha$. 
This shows that $a\in V_\alpha[\genG]$. 
\qedofLemma\qedskip\fi}

\begin{Lemma}\Label{p-Lg-RcA-0-1}
  For any $n\in\natnums\setminus1$, there is a $\Sigma_n$-formula
  $\varphi_n=\varphi_n(\overline{x},y)$ (for each finite sequence $\overline{x}$ of 
  variables\footnotemark) \st, if $\uniU$ is a transitive models of large enough fragment 
  of \ZFC, then  
  for any $\Sigma_n$ formula $\psi=\psi(\overline{x})$ there is $p\in\calH(\omega)$ \st\
  $\uniU\models\psi(\overline{a})$ if and only if $\uniU\models\varphi_n(\overline{a},p)$ 
  for all $a\in\uniU$. \ifextended\else\qed\fi
\end{Lemma}
\footnotetext{Similarly to the convention of some computer languages we consider here that 
  we have distinct $\varphi_n(\overline{x},y)$ for each length of sequence $\overline{x}$ 
  of variables.}
{\ifextended\extendedcolor\prf
  For $n=1$, $\varphi_1(\overline{x},y)$ can be chosen as a $\Sigma_1$-formula saying
  \begin{xitemize}
  \xitemA[x-2-0] 
    $\exists M\,(\,
    \begin{array}[t]{@{}l}
      M\mbox{ is transitive, } \overline{x}\in M,\,
      M\modelof{\mbox{large enough fragment of ZFC}}\\
      y\mbox{ is a code of a }\Sigma_1\mbox{-formula and }M\models\unmeralof{y}(\overline{x})
).
    \end{array}$
  \end{xitemize}

  If $\varphi_n(\overline{x},y)$ is defined for $\overline{x}$ of various lengths we can 
  define $\varphi_{n+1}(\overline{x})$ as $\exists x\neg\varphi_n(\overline{x},x,y)$. 
\qedofLemma\fi}

\begin{Lemma}
  \Label{p-Lg-RcA-0-2} {\ifextended\extendedcolor\wassertof{1} \fi} For 
  any $n^*\in\natnums$ there is $n>n^*$ \st, if\/ 
  $V_\alpha\prec_{\Sigma_n}\uniV$, then $V_\alpha[\genG]\prec_{\Sigma_{n^*}}\uniV[\genG]$ 
  for any $\poP\in V_\alpha$ and $(V,\poP)$-generic $\genG$. \ifextended\else\qed\fi
{\ifextended\extendedcolor\smallskip\\
\wassert{2} For a natural number $n$, there is $n'>n$ \st, for any $\alpha\in\On$, if\/ 
$V_\alpha[\genG]\prec_{\Sigma_{n'}}\uniV[\genG]$ for a \po\/ $\poP\in V_\alpha$ and
$(\uniV,\poP)$-generic $\genG$, then we have $V_\alpha\prec_{\Sigma_n}\uniV$. 
\fi}
\end{Lemma}
{\ifextended\extendedcolor\prf
\memo{{\tt recurrence-axioms.tex p-Lg-RcA-4-0 (1)}}
\assertof{1}: Suppose that $n>n^*$ is sufficiently large, $V_\lambda\prec\uniV$,
$\overline{a}\in V_\lambda[\genG]$, and 
$\varphi=\varphi(\overline{x})$ is  
a $\Sigma_{n^*}$-formula. 
There are $\poP$-names 
$\overline{\uta}\in V_\lambda$ \st\ $\overline{a}=\overline{\uta}[\genG]$.

If $V_\lambda[\genG]\models\varphi(\overline{a})$, there is $\condp\in\genG$ \st\
$V_\lambda\models\condp\forces{\poP}{\varphi(\overline{\uta})}$. By the choice of $n$ it 
follows that $\uniV\models\condp\forces{\poP}{\varphi(\overline{\uta})}$. 
Thus $\uniV[\genG]\models\varphi(\overline{a})$. 

The same argument also applies to $\neg\varphi$. \smallskip

\assertof{2}: We use the $\Lin$-formula $\Phi(x,y)$ of \Lemmaof{p-intro-1}. 
By assumption, there is $r\in V_\alpha[\genG]$ \st\ 
\begin{itemize}
\item[] $V_\alpha=\Phi(\cdot,r)^{V_\alpha[\genG]}
  =\Phi(\cdot,r)^\uniV\cap V_\alpha[\genG]\subseteq\Phi(\cdot,r)^\uniV$.
\end{itemize}
For any $\Sigma_n$-formula $\varphi(\overline{x})$ and
$\overline{a}\in\Phi(\cdot,r)^{V_\alpha[\genG]}$. Since $\varphi^{\Phi(\cdot,r)}$ is a
$\Sigma_{n'}$-formula (by the choice of $n'$), we have\vspace{-2ex}
\begin{xitemize}
\item[] $V_\alpha\models\varphi(\overline{a})$\ \ 
  $\Leftrightarrow$\ \  
  $V_\alpha[\genG]\models\varphi^{\Phi(\cdot,r)}(\overline{a})$\ \ 
  {$\obecause{\Leftrightarrow}{}{by assumption}$}\ \  
  $\uniV[\genG]\models\varphi^{\Phi(\cdot,r)}(\overline{a})$\\
  $\Leftrightarrow$\ \  
  $V\models\varphi(\overline{a})$. 
\end{xitemize}
This shows that $V_\alpha\prec_{\Sigma_n}\uniV$.
\qedofLemma\qedskip\fi}

\RcA\ and \RcAp\ 
are actually (almost) identical with (certain variations of) already well-known axioms 
and principles. 

For a class $\calP$ of \pos, an $\Lin$-formula $\varphi(\overline{a})$ with 
parameters 
$\overline{a}$  $(\in\uniV)$ is said to be a {\It$\calP$-button} if there is $\poP\in\calP$ 
\st, for any $\poP$-name $\utpoQ$ of \po\ with $\forces{\poP}{\utpoQ\in\calP}$, we have
$\forces{\poP\ast\utpoQ}{\varphi(\overline{a}\checked)}$.

If $\varphi(\overline{a})$ is a $\calP$-button then we call $\poP$ as above a {\It push of 
  the button $\varphi(\overline{a})$}. 

For a class $\poP$ of \pos\ and a set $A$ (of parameters), the {\It Maximality Principle 
for $\calP$ and $A$} ({\It$\MP(\calP,A)$}, for short) introduced in Hamkins \cite{hamkins} is 
the following assertion formulated in an axioms scheme in $\Lin$:

\begin{itemize}
\item[{\darkred$\MP(\calP,A)$}: ]\qquad\qquad For 
  any $\Lin$-formula $\varphi(\overline{x})$ and 
  $\overline{a}\in A$, if $\varphi(\overline{a})$ is a $\calP$-button then
  $\varphi(\overline{a})$ holds. 
\end{itemize}

\begin{Prop}{\rm(Barton, Caicedo, Fuchs, Hamkins, Reitz, and Schindler \cite{5a})}
  \Label{p-intro-0} Suppose that $\calP$  is an \jumpto{iterable}{iterable} class of \pos\ 
  and $A$ a set 
  (of parameters). \wassertof{1} $(\calP, A)$-\RcAp\ is equivalent to
  $\MP(\calP,A)$.\smallskip 

  \wassert{2} $(\calP, A)$-\RcA\ is equivalent to the following assertion: 
\begin{xitemize}
\xitem[x-a]
  For any $\Lin$-formula $\varphi(\overline{x})$ and
  $\overline{a}\in A$, if $\varphi(\overline{a})$ is a $\calP$-button then
  $\varphi(\overline{a})$ holds in a ground of\/ $\uniV$. 
\end{xitemize}
\end{Prop}
\prf \assertof{1}: Suppose first that $(\calP, A)$-\RcAp\ holds. We show that 
\memo{\MA\ $\rightarrow$ \MP}
$\MP(\calP,A)$ holds. Suppose that $\poP\in\calP$ is a push of the $\calP$-button
$\varphi(\overline{a})$. Let $\varphi'(\overline{x})$ be the formula expressing
\begin{equation}
\Label{x-b}
  \mbox{for any }\poQ\in\calP,\ 
  \forces{\poQ}{\varphi(\overline{x}\checked)}\mbox{ holds. }
\end{equation}
Then we have $\forces{\poP}{\varphi'(\overline{a}\checked)}$. By 
$(\calP, A)$-\RcAp, there is a $\calP$-ground $M$ of $\uniV$ \st\ $\overline{a}\in M$ and
$M\models\varphi'(\overline{a})$ holds. By the definition \xitemof{x-b} 
of $\varphi'$, it follows that $\uniV\models\varphi(\overline{a})$ holds. 

\memo{\MA\ $\rightarrow$ \MP}
Now suppose that $\MP(\calP,A)$ holds and $\poP\in\calP$ is \st\
$\forces{\poP}{\varphi(\overline{a})}$ 
for $\overline{a}\in A$.

Let $\varphi'$ be a formula claiming that
\begin{equation}
  \Label{x-c}
  \mbox{there is a }\calP\mbox{-ground }N\mbox{ \st\ }\overline{x}\in N\mbox{ and }
  N\models\varphi(\overline{x}). 
\end{equation}
Then $\varphi'(\overline{a})$ is a $\calP$-button and $\poP$ is its push.

\memo{\MA\ $\rightarrow$ \MP}
By $\MP(\calP,A)$, $\varphi'(\overline{a})$ holds in $\uniV$ and hence there is 
a $\calP$-ground $M$ 
of $\uniV$ \st\ $\overline{a}\in M$ and $M\models\varphi(\overline{a})$. This shows that 
$(\calP, A)$-\RcAp\ holds. \smallskip

\assertof{2}: can be proved similarly to \assertof{1}. 
Suppose first that $(\calP, A)$-\RcA\ holds. We show that 
\xitemof{x-a} 
holds. Suppose that $\poP\in\calP$ is a push of the $\calP$-button
$\varphi(\overline{a})$. Let $\varphi'(\overline{x})$ be the formula expressing
\begin{equation}
  \Label{x-d}
  \mbox{for any }\poQ\in\calP,\ 
  \forces{\poQ}{\varphi(\overline{x}\checked)}\mbox{ holds.}
\end{equation}
Then we have $\forces{\poP}{\varphi'(\overline{a}\checked}$. By 
$(\calP, A)$-\RcA, there is a ground $M$ of $\uniV$ \st\ $\overline{a}\in M$ and
$M\models\varphi'(\overline{a})$ holds. 
Since $\calP\ni\ssetof{\bbone}$, it follows that $M\models\varphi(\overline{a})$. 

Now suppose that \xitemof{x-a} 
holds and $\poP\in\calP$ is \st\ 
$\forces{\poP}{\varphi(\overline{a}\checked)}$ 
for $\overline{a}\in A$.

Let $\varphi''$ be a formula asserting that
\begin{equation}
  \Label{x-e}
  \mbox{there is a }\calP\mbox{-ground }N\mbox{ \st\ }\overline{x}\in N\mbox{ and }
  N\models\varphi(\overline{x}). 
\end{equation}
Then $\varphi''(\overline{a})$ is a $\calP$-button and $\poP$ is its push.
Thus, By 
\xitemof{x-a}, $\varphi''(\overline{a})$ holds in a ground $M$ of $\uniV$ 
with $\overline{a}\in M$. By the definition 
\xitemof{x-e} 
of $\varphi''$, there is a $\calP$-ground $N$ of $M$ \st\
$\overline{a}\in N$ and $N\models\varphi(\overline{a})$. Since $N$ is also a ground of
$\uniV$, this shows that 
$(\calP, A)$-\RcA\ holds. 
\qedoftheorem\qedskip

Recurrence Axioms are also related to the Inner Model Hypothesis introduced by Sy 
Friedman in \cite{friedman-sy}. {\It The Inner Model Hypothesis} ({\It\IMH}) is the 
following assertion formulated in the language of second-order set theory (e.g.\ 
in the context of von Neumann-Bernays-Gödel set theory): 
\begin{itemize}
\item[{\It\IMH} :]\qquad For any statement $\varphi$ without parameters, if $\varphi$ holds 
  in an  
  inner model of an inner extension of $\uniV$ then $\varphi$ holds in an inner model of
  $\uniV$. 
\end{itemize}
Here we say a (not necessarily first-order definable) transitive class $M$ an {\It inner 
  model}\/  
of $\uniV$ if $M$ is a model of \ZF\ and $\On^M=\On^\uniV$. In the perspective from 
such $M$, we call $\uniV$ an {\It inner extension} of $M$. 

We call a set-forcing version of this principle {\It Inner Ground Hypothesis} 
({\It\IGH}):  

For a (definable) class $\calP$ of \pos\ and a set $A$ (of parameters), 
\begin{itemize}
\item[{\It$\IGH(\calP,A)$} :]\qquad\qquad\ \ For any $\Lin$-formula
  $\varphi=\varphi(\overline{x})$ and 
  $\overline{a}\in A$, if $\poP\in\calP$ forces ``there is a ground $M$ with
  $\overline{a}\in M$ satisfying $\varphi(\overline{a})$'', then there is a ground $\uniW$ of
  $\uniV$ \st\ $\overline{a}\in\uniW$ and $\uniW\models\varphi(\overline{a})$. 
\end{itemize}

\begin{Prop}{\rm(Barton, Caicedo, Fuchs, Hamkins, Reitz, and Schindler \cite{5a})}
  \Label{p-intro-2} For a class $\calP$ of \pos\ with $\ssetof{\bbone}\in\calP$ and a set $A$ (of 
  parameters), $(\calP,A)$-\RcA\ holds if and only if\/ $\IGH(\calP,A)$ holds.
\end{Prop}
\prf Suppose that $(\calP,A)$-\RcA\ holds. Let $\varphi=\varphi(\overline{x})$ be 
an $\Lin$-formula, $\overline{a}\in A$, and $\poP\in\calP$  be \st\
$\forces{\poP}{\varphi(\overline{a}\checked)\xmbox{ holds in a ground}}$. 

Let $\varphi'(\overline{x})$ be the $\Lin$-formula asserting that $\varphi(\overline{x})$ 
holds in a ground. Then $\forces{\poP}{\varphi'(\overline{a}\checked)}$. 
By $(\calP,A)$-\RcA, it follows that there is a ground $\uniW$ of $\uniV$ \st\
$\uniW\models\varphi'(\overline{a}\checked)$. Since a ground of a ground is a ground, we 
conclude that there is a ground $\uniW_0$ of $\uniV$ \st\ $\overline{a}\in M_0$ and
$\uniW_0\models\varphi(\overline{a})$. This shows that $\IGH(\calP,A)$ holds. 

Suppose now that $\IGH(\calP,A)$ holds. Assume that
$\forces{\poP}{\varphi(\overline{a}\checked)}$ for an $\Lin$-formula
$\varphi=\varphi(\overline{x})$, $\overline{a}\in A$, and $\poP\in\calP$. Then
$\forces{\poP}{\varphi(\overline{a}\checked)\xmbox{ holds in a }
\calP\xmbox{-ground (of the universe)}}$ since $\forces{\poP}{\ssetof{\bbone}\in\calP}$. 
Thus, by $\IGH(\calP,A)$, there is a ground $\uniW$ of $\uniV$ \st\
$\uniW\models\varphi(\overline{a})$. 
\qedofProp
\qedskip

\memo{\normalsize!!!}
{\extendedcolor$(\calP,A)$}-\RcAp\ ($\Leftrightarrow$\ $\MP(\calP,A)$ for an iterable $\calP$) can be also 
characterized in terms of a strengthening of Inner Ground Hypothesis:
For a (definable) class $\calP$ of \pos\ and a set $A$ (of parameters), 
\begin{itemize}
\item[{\It$\IGH^+(\calP,A)$} :]\qquad\qquad\quad For any $\Lin$-formula
  $\varphi=\varphi(\overline{x})$ and 
  $\overline{a}\in A$ if $\poP\in\calP$ forces ``there is a $\calP$-ground $M$ with
  $\overline{a}\in M$ satisfying $\varphi(\overline{a})$'', then there is a $\calP$-ground
  $\uniW$ of 
  $\uniV$ \st\ $\overline{a}\in\uniW$ and $\uniW\models\varphi(\overline{a})$. 
\end{itemize}

The following proposition can be proved similarly to \Propof{p-intro-2}.

\begin{Prop}
  \Label{p-intro-3} For an iterable class $\calP$ of \pos\ and a set $A$ (of parameters),
  $(\calP,A)$-\RcAp\ holds if and only if $\IGH^+(\calP,A)$ holds. 
  \ifextended\else\qed\fi
\end{Prop}
{\ifextended\extendedcolor
\prf Suppose that $(\calP,A)$-\RcAp\ holds and assume that $\varphi=\varphi(\overline{x})$ 
is an $\Lin$-formula, $\overline{a}\in A$, and $\poP\in\calP$ is \st\ 
\begin{equation*}
  \forces{\poP}{\varphi(\overline{a}\checked)\mbox{ holds in a }\calP\mbox{-ground }M
  \mbox{ with }\overline{a}\in M}
\end{equation*}
Let $\varphi'(\overline{a})$ be the formula expressing ``$\varphi(\overline{x})$ holds in a
$\calP$-ground $M$ with $\overline{a}\in M$''. Then $\poP$ is a push of the switch
$\varphi'(\overline{a})$. Thus, by \Propof{p-intro-0},\,\assertof{1},
$\varphi'(\overline{a})$ holds in $\uniV$. By definition of $\varphi'$, there is 
a $\calP$-ground $\uniW$ of $\uniV$ \st\ $\overline{a}\in\uniW$ and
$\uniW\models\varphi(\overline{a})$. 
This shows that
$\IGH^+(\calP,A)$ holds. 

Suppose now that $\IGH^+(\calP,A)$ holds, and assume that $\varphi=\varphi(\overline{x})$ 
is an $\Lin$-formula, $\overline{a}\in A$ 
and $\forces{\poP}{\varphi(\overline{a}\checked)}$ then  
(since $\ssetof{\bbone}\in\calP$)
$\forces{\poP}{\varphi(\overline{a}\checked)
  \xmbox{ holds in a }\calP\xmbox{-ground }M\xmbox{ with }\overline{a}\in M}$.
By $\IGH^+(\calP,A)$, it follows that there is  $\calP$-ground $\uniW_0$ of $\calP$-ground of
$\uniV$ \st\  
$\overline{a}\in\uniW_0$ and $\uniW_0\models\varphi(\overline{a})$. Since $\calP$ is 
iterable,  
$\uniW_0$ is a $\calP$-ground of $\uniV$. 
This shows 
that $(\calP,A)$-\RcAp\ holds. 
\qedofProp
\qedskip\fi}

In spite of these characterizations and near characterizations, we want to keep the 
Recurrence Axioms as autarchic axioms. On of the reasons is that we want to retain the 
narration that is is a reflection principle with reflection from the set-theoretic 
multiverse to the geology; another reason is that we have the following 
monotonicity which does not hold e.g.\ for Maximality Principles. 

\memo{!!! lemma -> Lemma}
\begin{Lemma}{\rm (Monotonicity of Recurrence Axioms)}\Label{p-5}
  For classes of \pos\ $\calP$, $\calP'$ and sets $A$, $A'$ of parameters, 
  if $\calP\subseteq \calP'$ and $A\subseteq A'$, then we have 
  \begin{xitemize}
  \item[] $(\calP', A')$-\RcA\ \ $\Rightarrow$\ \ $(\calP,A)$-\RcA. \qed
  \end{xitemize}
\end{Lemma}

If we decide that the Recurrence Axioms are desirable extensions of the axioms of 
\ZFC, then we should adopt the  maximal instance of these axioms. (i.e.\ the one with 
maximal strength among the 
instances consistent with \ZFC) 
By \Lemmaof{p-5}, this means we should try to take the instance of Recurrence Axioms  
with the maximal $\calP$ and $A$ (\wrt\ inclusion) among the consistent ones. 

\theoremof{p-Lg-RcA-1} \memox{\large !!!} in the next section suggests that the following two 
as candidates of such maximal instances:
\begin{itemize}
\item[\wassertof{$E$}]\qquad\ $\ZFC$ $+$ $(\calP,\calH(\kappa_\refl))$-RcA\ for the class 
  $\calP$ of all  
  stationary preserving \pos.
\end{itemize}\begin{itemize}
\item[\wassertof{$Z$}]\qquad\ $\ZFC$ $+$ $(\calQ,\calH(\continuum))$-RcA\ for the class
  $\calQ$ of all \pos.  
\end{itemize}

The consistency of \assertof{$Z$} follows from the consistency of \ZFC\ $+$ ``there are 
stationarily many inaccessible cardinals'' (\cite{hamkins}). The consistency of 
\assertof{$E$} follows from  \Lemmaof{p-Lg-RcA-2}, 
\theoremof{p-Lg-RcA-4},\,\assertof{$B'$}, and   
\theoremof{p-Lg-RcA-5}.\memox{\large!!!}

The maximality of \assertof{$E$} and \assertof{$Z$} follows from 
\Lemmaof{p-Lg-RcA-1},\,\assertof{2'} and \assertof{5'} respectively.\memox{\large!!!}

By \Lemmaof{p-Lg-RcA-1},\,\assertof{4} and \assertof{5}, \assertof{$E$} implies
$\continuum=\aleph_2$, and \assertof{$Z$} implies \CH. In  
particular, these two extensions of \ZFC\ are not compatible. However, as we are going to 
discuss in \sectionof{sec:5}, we can combine \assertof{$E$} with a reasonable weakening of 
\assertof{$Z$}. 

\begin{itemize}\Label{z+}
\item[\wassertof{$Z^+$}]\qquad\ $\ZFC$ $+$ $(\calP,\calH(\kappa_\refl))$-\RcAp\ $+$ 
$(\calQ,{\calH(\omega_1)}^{\overline{\uniW}})$-\RcAp\ where $\calP$ is the class of all 
  proper \pos, $\calQ$ the class of all \pos, and $\overline{\uniW}$ the 
  bedrock\footnotemark\ which is also assumed here to exist.
\end{itemize}
\footnotetext{For the definition of the bedrock see \sectionof{sec:3}.}

\assertof{$Z^+$} implies $\continuum=\aleph_2$ (see \Lemmaof{p-Lg-RcA-1}, \assertof{4} in 
the next section).  
This is what I meant when I wrote ``the maximal 
setting of Recurrence Axioms points to the universe with the continuum of size $\aleph_2$'' 
in the introduction. We shall further discuss about \assertof{$Z^+$} in \sectionof{sec:4}. 


\section{Restricted Recurrence Axioms}
\Label{sec:2-0} The following restricted forms of Recurrence Axioms are enough to 
decide many interesting aspects including the cardinal arithmetic around the continuum. 

For an iterable class $\calP$ 
of \pos, a set $A$ (of parameters), and a set $\Gamma$ of $\Lin$-
formulas, {\It$\calP$-Recurrence Axiom for formulas in $\Gamma$ with parameters from $A$} 
({\It$(\calP,A)_\Gamma$-\RcA}, for short) is the following assertion expressed as an axiom 
scheme in $\Lin$:
\begin{itemize}
\item[{\It$(\calP,A)_\Gamma$-\RcA}:]\qquad\qquad\quad For any $\varphi(\overline{x})\in\Gamma$ and $\overline{a}\in A$, if
  $\forces{\poP}{\varphi(\overline{a}\checked)}$, then there is a ground 
  $\uniW$ of $\uniV$ \st\ $\overline{a}\in \uniW$ and $\uniW\models\varphi(\overline{a})$. 
\end{itemize}
{\ifextended\else\It\fi$(\calP,A)_\Gamma$-\RcAp} corresponding to $(\calP,A)$-\RcAp\ is 
defined similarly.  

\ifextended
\begin{itemize}
\item[{\It$(\calP,A)_\Gamma$-\RcAp}:]\qquad\qquad\quad For any $\varphi(\overline{x})\in\Gamma$ and $\overline{a}\in A$, if
  $\forces{\poP}{\varphi(\overline{a}\checked)}$, then there is a\\ $\calP$-ground 
  $\uniW$ of $\uniV$ \st\ $\overline{a}\in \uniW$ and $\uniW\models\varphi(\overline{a})$. 
\end{itemize}
\fi
\begin{Lemma}{\rm(Fuchino and Usuba \cite{laver-gen-maximum})}
  \Label{p-Lg-RcA-1} Assume that $\calP$ is an \jumpto{iterable}{iterable} class of \pos. \wassertof{1} If 
  $\calP$ contains a \po\ which adds a real (over 
  the universe), then 
  $(\calP,\calH(\kappa_\refl))_{\Sigma_1}$-\RcA\ implies $\neg\CH$.\smallskip

  \wassert{2} Suppose that $\calP$ contains a 
  \po\ which 
  forces ${\aleph_2}^\uniV$ to be equinumerous with ${\aleph_1}^\uniV$. Then 
  $(\calP,\calH(2^{\aleph_0}))_{\Sigma_1}$-\RcA\ implies $\continuum\leq\aleph_2$.
  \smallskip

  \wassert{2'} If $\calP$ contains a \pos\ which forces ${\aleph_2}^\uniV$ to be 
  equinumerous with 
  ${\aleph_1}^\uniV$, then $(\calP,\calH((\aleph_2)^+))_{\Sigma_1}$-\RcA\ does not hold. 
  \smallskip

  \wassert{3} If $(\calP,\calH(\kappa_\refl))_{\Sigma_1}$-\RcA\ holds then all 
  $\poP\in\calP$ preserve $\aleph_1$ and they are also stationary preserving.\smallskip

  \wassert{4} If $\calP$ contains a \po\ which 
  adds a real as well as a \po\ which collapses 
  ${\aleph_2}^\uniV$, then 
  $(\calP,\calH(\kappa_\refl))_{\Sigma_1}$-\RcA\ implies $\continuum=\aleph_2$.
  \smallskip

  \wassert{5} If $\calP$ contains a 
  \po\ which 
  collapses ${\aleph_1}^\uniV$, then $(\calP,\calH(2^{\aleph_0}))_{\Sigma_1}$-\RcA\ implies 
  \CH.\smallskip 

  \wassert{5'} If $\calP$ contains a \po\ which collapses ${\aleph_1}^\uniV$ then 
  $(\calP,\calH((\continuum)^+))_{\Sigma_1}$-\RcA\ does not hold. 

  {\ifextended\extendedcolor\wassert{6} Suppose that all $\poP\in\calP$ preserve cardinals 
    and $\calP$ contains \pos\  
  adding at least $\kappa$ many reals for eachi $\kappa\in\Card$. Then
  $(\calP,\emptyset)_{\Sigma_2}$-\RcAp\ implies that $\continuum$ is very large.

  \wassert{6'} Suppose that $\calP$ is as in \assertof{6}. Then 
  $(\calP,\calH(\continuum))_{\Sigma_2}$-\RcAp\ implies that $\continuum$ is a limit 
  cardinal. Thus if $\continuum$ is regular in addition, then $\continuum$ is weakly 
  incaccessible.\fi} 
\end{Lemma}
\ifextended{\extendedcolor
\prf \assertof{1}: Assume that $\calP$ is an iterable class of \pos\ containing a 
\po\ $\poP$ adding a real and $(\calP,\calH(\kappa_\refl))_{\Sigma_1}$-\RcA\ holds. If \CH\ 
holds, then $\psof{\omega}^\uniV\in\calH(\kappa_\refl)$. Hence 
\begin{xitemize}
\xitemA[x-3]
  $\mbox{``\,}\exists x\,(x\subseteq\omega\land x\not\in\psof{\omega}^\uniV)\mbox{''}$ 
\end{xitemize}
is 
a $\Sigma_1$-formula with parameters from $\calH(\kappa_\refl)$ and $\poP$ forces the 
formula in the forcing language corresponding to this 
formula:
``$\exists x\,(s\subseteq\check{\omega}\,\land x\not\in (\psof{\omega}^\uniV)\checked)$''.  

By $(\calP,\calH(\kappa_\refl))_{\Sigma_1}$-\RcA, the formula \xitemAof{x-3} must hold in 
a ground. This is a contradiction. \smallskip

\assertof{2}: Assume that $(\calP,\calH(2^{\aleph_0})_{\Sigma_1})$-\RcA\ holds and
$\poP\in\calP$ 
forces ${\aleph_2}^\uniV$ to be equinumerous with ${\aleph_1}^\uniV$. 
If $\continuum>\aleph_2$ then 
${\aleph_1}^\uniV$, ${\aleph_2}^\uniV\in\calH(2^{\aleph_0})$. Letting 
$\psi(x,y)$ a $\Sigma_1$-formula saying 
``$\exists f\,(f\mbox{ is a surjection from }x\mbox{ to }y)$'', 
we have 
$\forces{\poP}{\psi(({\aleph_1}^\uniV)\checked,
({\aleph_2}^\uniV)\checked)}$. 

By $(\calP,\calH(2^{\aleph_0}))_{\Sigma_1}$-\RcA, the formula 
$\psi({\aleph_1}^\uniV,{\aleph_2}^\uniV)$
must hold in a ground. This is a contradiction. \smallskip

\assertof{2'}: Assume that $\poP\in\calP$ is \st\
$\forces{\poP}{\cardof{({\aleph_2}^\uniV)\checked}=\cardof{({\aleph_1}^\uniV)\checked}}$, 
and 
$(\calP,\calH({\aleph_2}^+))_{\Sigma_1}$-\RcA\ holds. Then, since $\aleph_1$,
$\aleph_2\in\calH({\aleph_2}^+)$ and ``$\cardof{x}=\cardof{y}$'' is $\Sigma_1$, there is a 
ground $\uniW$ of $\uniV$ \st\
$\uniW\models\cardof{{\aleph_2}^\uniV}=\cardof{{\aleph_1}^\uniV}$. This is a 
contradiction.\smallskip

\assertof{3}: Suppose that $\poP\in\calP$ is \st\
$\forces{\poP}{{\aleph_1}^\uniV\mbox{ is countable}}$. Note that
$\omega, \aleph_1\in\calH(\kappa_\refl)$. 
By
$(\calP,\calH(\kappa_\refl))_{\Sigma_1}$-\RcA, it follows that there is a ground $\uniW$ of 
$\uniV$ \st\ $\uniW\modelof{{\aleph_1}^\uniV\mbox{ is countable}}$. This is a contradiction.

Suppose now that $S\subseteq\omega_1$ is stationary and $\poP\in\calP$ destroys the 
stationarity of $S$. Note that 
$\omega_1$, $S\in\calH(\aleph_2)$. Let $\varphi=\varphi(y,z)$ be the $\Sigma_1$-formula
\begin{itemize}
\item[] $\exists x\,(x\mbox{ is a club subset of the ordinal }y\mbox{ and }
  z\cap x=\emptyset)$. 
\end{itemize}
Then we have $\forces{\poP}{\varphi(\omega_1, S)}$. By
$(\calP,\calH(\kappa_\refl))_{\Sigma_1}$-\RcA, it follows that there is a ground
$\uniW\subseteq\uniV$ \st\ $S\in\uniW$ and $\uniW\models\varphi(\omega_1,S)$. This is a 
contradiction to the stationarity of $S$. 

\smallskip

\assertof{4}: follows from \assertof{1}, \assertof{2} and \assertof{3}. 
\smallskip

\assertof{5}: 
If $\aleph_1<\continuum$, 
then ${\aleph_1}^\uniV\in\calH(2^{\aleph_0})$.

Let $\poP\in\calP$ be a \po\ collapsing ${\aleph_1}^\uniV$. That is, 
$\forces{\poP}{{\aleph_1}^\uniV\mbox{ is countable}}$. Since ``$\cdots$ is countable'' is
$\Sigma_1$, there is a ground $M$ \st\ $M\modelof{{\aleph_1}^\uniV\mbox{ is countable}}$ by
$(\calP, \calH(2^{\aleph_0}))_{\Sigma_1}$-\RcA.  
This is a contradiction. \smallskip

\assertof{5'}: Assume that $\poP\in\calP$ is \st\
$\forces{\poP}{{\aleph_1}^\uniV\mbox{ is countable}}$, 
and $(\calP,\calH((\continuum)^+))$-\RcA\ holds. Since $\aleph_1\in\calH((\continuum)^+)$, 
it follows that there is a ground $\uniW$ of $\uniV$ \st\
$\uniW\models{\aleph_1}^\uniV\mbox{ is countable}$. This is a contradiction.\smallskip

\assertof{6}: To prove e.g.\ that $\continuum>\aleph_\omega$, let
  $\poP\in\calP$ be \st\  
  $\forces{\poP}{\continuum>\aleph_\omega}$. Then by
  $(\calP,\emptyset)_{\Sigma_2}{\mbox{-}}$\RcAp, there is a $\calP$-ground $\uniW$ 
  of $\uniV$ \st\ $\uniW\models\continuum>\aleph_\omega$. Since $\uniV$ is $\calP$-gen.\ 
  extension of $\uniW$ and $\calP$ preservs cardinals, it follows that
  $\uniV\models\continuum>\aleph_\omega$.\smallskip

  \assertof{6'}: Suppose $\mu<\continuum$. Then $\mu\in\calH(\continuum)$. There is $\poP\in\calP$ 
  \st\ $\forces{\poP}{\continuum>\mu^+}$. By
  $(\calP,\calH(\continuum))_{\Sigma_2}\mbox{-}$\RcAp, it follows that there is a $\calP$-ground 
$\uniW$ of $\uniV$ which satisfies this statement. Since $\calP$ preserves cardinals it 
  follows that $\uniV\models\continuum>\mu^+$. 
\qedofLemma\qedskip}\else\prf All assertions are proved similarly. We see here only the 
proof of \assertof{1}: Assume that $\calP$ is an iterable class of \pos\ containing a 
\po\ $\poP$ adding a real and $(\calP,\calH(\kappa_\refl))_{\Sigma_1}$-\RcA\ holds. If \CH\ 
holds, then $\psof{\omega}^\uniV\in\calH(\kappa_\refl)$. Hence 
\begin{equation}\Label{x-3}
  \mbox{``\,}\exists x\,(x\subseteq\omega\land x\not\in\psof{\omega}^\uniV)\mbox{''} 
\end{equation}
is 
a $\Sigma_1$-formula with parameters from $\calH(\kappa_\refl)$ and $\poP$ forces the 
formula in the forcing language corresponding to this 
formula:
``$\exists x\,(s\subseteq\check{\omega}\,\land x\not\in (\psof{\omega}^\uniV)\checked)$''.  

By $(\calP,\calH(\kappa_\refl))_{\Sigma_1}$-\RcA, the formula \equationof{x-3} must hold in 
a ground. This is a contradiction.

\qedofLemma\qedskip
\fi

{\ifextended\extendedcolor
  \begin{CorA}\extendedcolor
    \Label{p-5-0}\wassertof{$E$} $(\mbox{stationary preserving},\calH(\kappa_\refl))$-\RcA\ 
    is maximal among Recurrence Axioms with similar pair of parameters, and it implies
    $\continuum=\aleph_2$. \smallskip

    \wassert{$Z$} $(\mbox{all \pos},\calH(\continuum))$-\RcA\ is maximal among Recurrence 
    Axioms with similar pair of parameters, and it implies $\CH$.
  \end{CorA}
  \prf \assertof{$E$}: By \Lemmaof{p-Lg-RcA-1},\,\assertof{3}, ``stationary preserving'' 
  cannot be replaced by a larger class of \pos. By \Lemmaof{p-Lg-RcA-1},\,\assertof{4},
  $(\mbox{stationary preserving},\calH(\continuum))$-\RcA\ implies $\continuum=\aleph_2$.
  Thus $\calH(\kappa_\refl)=\calH(\aleph_2)$ in this case and, by 
  \Lemmaof{p-Lg-RcA-1},\,\assertof{2'}, this cannot be replaced 
  by $\calH(\calH(\aleph_3))$. \smallskip

  \assertof{$Z$}: \CH\ holds by \Lemmaof{p-Lg-RcA-1},\,\assertof{5}, and, by By 
  \Lemmaof{p-Lg-RcA-1},\,\assertof{5'}, $\calH(\continuum)$ cannot be replaced by
  $\calH((\continuum)^+)$. \qedofCorA\qedskip
\fi}

Laver-genericity implies the plus version of Recurrence Axiom ($\Leftrightarrow$ 
Maximality Principle) restricted to $\Sigma_2$.

\begin{theorem}
  \Label{p-Lg-RcA-0} Suppose that $\kappa$ is \jumpto{LGHH}{tightly $\calP$-Laver-generically ultrahuge} 
  for an 
  \jumpto{iterable}{iterable} class $\calP$ of \pos
  . Then
  $(\calP,\calH(\kappa))_{\Sigma_2}$-\RcAp\ holds. 
\end{theorem}
\prf Assume that $\kappa$ is tightly $\calP$-Laver generically ultrahuge for an 
iterable class $\calP$ of \pos. 

Suppose that $\varphi=\varphi(\overline{x})$ is $\Sigma_2$ formula (in $\Lin$),
$\overline{a}\in\calH(\kappa)$, and  
$\poP\in\calP$ is \st\ 
\begin{equation}
\Label{x-Lg-RcA-a} 
  \uniV\models\forces{\poP}{\varphi(\overline{a}\checked)}. 
\end{equation}

Let $\lambda>\kappa$ be 
\st\ $\poP\in\uniV_\lambda$ and 
\begin{equation}
\Label{x-Lg-RcA-0} V_\lambda\prec^{}_{\Sigma_n}\uniV\mbox{ for a sufficiently large }n. 
\end{equation}
In particular, we may assume that we have chosen the $n$ above so that a sufficiently large 
fragment of \ZFC\ holds in 
$V_\lambda$ (``sufficiently large'' means here, in particular, in terms of 
\Lemmaof{p-Lg-RcA-0-0} and that the argument at the end of this proof is possible).  

Let $\utpoQ$ be a $\poP$-name \st\ $\forces{\poP}{\utpoQ\in\calP}$, and for
$(\uniV,\poP\ast\utpoQ)$-generic $\genH$, there are $j$, $M\subseteq\uniV[\genH]$ with 
\begin{equation}
  \Label{x-Lg-RcA-1} 
  \Elembed{j}{\uniV}{M}{\kappa}, 
\end{equation}
\begin{equation}
  \Label{x-Lg-RcA-1-0} 
  j(\kappa)>\lambda, 
\end{equation}
\begin{equation}
  \Label{x-Lg-RcA-1-1} 
  \poP\ast\utpoQ,\ \poP,\ \genH,\ {V_{j(\lambda)}}^{\uniV[\genH]}\in M,\mbox{ and }
\end{equation}
\begin{equation}
  \Label{x-Lg-RcA-1-2}
  \cardof{\poP\ast\utpoQ}\leq j(\kappa). 
\end{equation}
By \xitemof{x-Lg-RcA-1-2}, we may assume that the underlying set 
of $\poP\ast\utpoQ$ is $j(\kappa)$ and $\poP\ast\utpoQ\in {V_{j(\lambda)}}^\uniV$. 

Let $\genG:=\genH\cap\poP$. 
Note that $\genG\in M$ by \xitemof{x-Lg-RcA-1-1} and 
we have 
\begin{equation}
\Label{x-Lg-RcA-2} 
  {V_{j(\lambda)}}^M\ubecause{=}{}{by \xitemof{x-Lg-RcA-1-1}}
  {V_{j(\lambda)}}^{\uniV[\genH]}
  \obecause{=}{1.44ex}{\qquad\qquad\qquad
    \vbox{\hbox{Since ${V_{j(\lambda)}}^M$ ($={V_{j(\lambda)}^{\uniV[\genH]}}$) satisfies a 
    sufficiently large fragment of 
    \ZFC\vspace{-0.18ex}}\hbox{by elementarity of $j$, and hence the equality follows by 
        \Lemmaof{p-Lg-RcA-0-0}}}} 
  {V_{j(\lambda)}}^\uniV[\genH]. \qquad\qquad
\end{equation}
Thus, by \xitemof{x-Lg-RcA-1-1} and by the definability of grounds, we have 
${V_{j(\lambda)}}^\uniV\in M$ and ${V_{j(\lambda)}}^\uniV[\genG]\in M$. 
\begin{Claim}
  \Label{cl-Lg-RcA-0} ${V_{j(\lambda)}}^\uniV[\genG]\models\varphi(\overline{a})$. 
\end{Claim}
\noindent
\prfofClaim
By \Lemmaof{p-Lg-RcA-0-0}, ${V_\lambda}^{\uniV}[\genG]={V_\lambda}^{\uniV[\genG]}$, and
${V_{j(\lambda)}}^{\uniV}[\genG]={V_{j(\lambda)}}^{\uniV[\genG]}$. 
\memo{\normalsize!!!}
By \xitemof{x-Lg-RcA-0}, both ${V_\lambda}^\uniV[\genG]$ and $V_{j(\lambda)}^\uniV[\genG]$ 
satisfy large enough fragment of \ZFC. Thus 
\begin{equation}
  \Label{x-Lg-RcA-2-0} 
  {V_\lambda}^\uniV[\genG]\prec_{\Sigma_1}{V_{j(\lambda)}}^\uniV[\genG]. 
\end{equation}
By \xitemof{x-Lg-RcA-a} and \xitemof{x-Lg-RcA-0}, we have
${V_\lambda}^\uniV[\genG]\models\varphi(\overline{a})$. By \xitemof{x-Lg-RcA-2-0} and since 
$\varphi$ is $\Sigma_2$, it follows that
${V_{j(\lambda)}}^\uniV[\genG]\models\varphi(\overline{a})$.  
\qedofClaim\qedskip

Thus we have 
\begin{equation}
  \Label{x-Lg-RcA-3} 
  M\modelof{\mbox{there is a }
  \calP\mbox{-ground }N\mbox{ of }V_{j(\lambda)}\mbox{ with }N\models\varphi(\overline{a})}.
\end{equation}

By the elementarity \xitemof{x-Lg-RcA-1}, it follows that 
\begin{equation}
  \Label{x-Lg-RcA-4} 
  \uniV\modelof{\mbox{there is a }
  \calP\mbox{-ground }N\mbox{ of }V_{\lambda}\mbox{ with }N\models\varphi(\overline{a})}.
\end{equation}

Now by \xitemof{x-Lg-RcA-0}, it follows that there is a $\calP$-ground $\uniW$ of $\uniV$ \st\
$\uniW\models\varphi(\overline{a})$. \qedoftheorem\qedskip

{\ifextended\extendedcolor 
  \begin{CorA}\extendedcolor
    Suppose that $\kappa=\kappa_\refl$ is \jumpto{LGHH}{tightly $\calP$-Laver-generically 
      ultrahuge}  
  for an 
  \jumpto{iterable}{iterable} class $\calP$ of \pos
  . Then elements of $\calP$ are stationary preserving. 
  \end{CorA}
  \prf By \lemmaof{p-Lg-RcA-1} and \theoremof{p-Lg-RcA-0}.\qedofCorA
\fi}

\section{Recurrence, Laver-generic large cardinal, and beyond}
\Label{sec:3}

Laver-genericity does not imply full Recurrence Axiom (\theoremof{p-10} and \Corof{p-10-0}
below\memox{\large!!!}).  

\Hypertarget{phistar}{}
For an $\Lin$-formula $\psi=\psi(\overline{x})$, a (large) cardinal $\kappa$ 
is {\It $\psi$-absolute} if the formula $\psi$ is absolute between $V_\kappa$ and $\uniV$ 
(i.e.\
$(\forall\overline{x}\in V_y)(\psi^{V_y}(\overline{x})\leftrightarrow\psi(\overline{x}))$ 
holds for $y=\kappa$).
\begin{Lemma}
  \Label{p-6} For any $n\in\natnums$, there is an $\Lin$-formula $\psi^*_n$ \st, for any 
  inaccessible $\kappa$, $\kappa$ is $\psi^*_n$-absolute if and only if 
  \begin{xitemize}
  \xitem[x-4] for any ground $\uniW$ of $\uniV$ \st\ $\uniV=\uniW[\genG]$ for a \po\
    $\poP\in {V_\kappa}^\uniW$ and $(\uniW,\poP)$-generic $\genG$,\footnotemark\ 
    we have that all $\Sigma_n$-formulas are absolute between ${V_\kappa}^\uniW$ and $\uniW$.
  \end{xitemize}
\end{Lemma}
\footnotetext{Note that this includes the case that $\poP=\ssetof{\bbone}$ 
  and $\uniV=\uniW$.}%
\prf By \theoremof{p-intro-1} and \Lemmaof{p-Lg-RcA-0-1}.
\qedofLemma
\begin{Lemma}
  \Label{p-7} Let $\psi^*_n$ be as in \Lemmaof{p-6}. $\psi^*_n$-absolute inaccessible 
  cardinals are not resurrectable. That is, if a cardinal $\lambda$ satisfies
  \begin{equation}
    \Label{x-5} \forces{\poP}{\check{\lambda}\mbox{ is }\psi^*_n\mbox{-absolute 
        inaccessible}} 
  \end{equation}
  for some \po\ $\poP$, then $\lambda$ is really $\psi^*_n$-absolute inaccessible. 
\end{Lemma}
\prf \equationof{x-5} implies that $\lambda$ is inaccessible. \memox{\normalsize!!!}
By the definition \xitemof{x-4} of $\psi^*_n$, if $\psi^*_n$ is absolute between 
${V_\lambda}^{\uniV[\genG]}$ and $\uniV[\genG]$ for some $(\uniV,\poP)$-generic $\genG$ 
then, it is absolute between ${V_\lambda}^\uniV$ and $\uniV$. 
\qedofLemma

\begin{Lemma}
  \Label{p-8} Suppose that there are stationarily many inaccessible cardinals.\footnotemark
  Then, for each $n\in\natnums$, there are stationarily many $\phi^*_n$-absolute inaccessible 
  cardinals. 
\end{Lemma}
\footnotetext{``There are stationarily many inaccessible cardinals'' is the statement 
  formalizable in an axiom scheme claiming, for 
  each $\Lin$-formula $\varphi=\varphi(x)$, that ``if $\varphi(x)$ defines a club subclass of
  $\On$  then there is an inaccessible $\mu$ with $\varphi(\mu)$''.  }
\prf
For $n\in\natnums$ let $n^+\geq n$ be \st\ $\psi^*_n$ is $\Sigma_{n^+}$. For any club
$C\subseteq\On$, $C\cap C^{(n^+)}=\setof{\alpha\in C}{V_\alpha\prec_{\Sigma_{n^+}\uniV}}$ is 
a club in $\On$ (Lévy-Montague Reflection Theorem), there is an inaccessible cardinal 
$\mu\in C\cap C^{(n^+)}$. By the choice of $n^+$, such $\mu$ is a $\psi^*_n$-absolute 
inaccessible cardinal.
\qedofLemma
\begin{theorem}
  \Label{p-9} Suppose that $(\calP,\emptyset)$-\RcA\ holds, where $\calP$ is a class of 
  \pos\ \st\ either \assertof{a} $\calP$ contains \pos\ collapsing arbitrary large 
  cardinals to a small cardinality (less than the first inaccessible if there are 
  inaccessibles at all), or \assertof{b} $\calP$ contains \pos\ adding arbitrarily many 
  reals.

  If there is a $\psi^*_n$-absolute inaccessible cardinal for some $n\in\natnums$, then 
  there are cofinally many $\psi^*_n$-absolute inaccessible cardinals.
\end{theorem}
\prf Assume that $(\calP,\emptyset)$-\RcA\ holds for $\calP$ as above and there is a 
cardinal $\lambda$ \st\ there are some $\psi^*_n$-absolute inaccessible cardinals but all 
of them are below $\lambda$.

Let $\poP$ be a \po\ which either collapses $\lambda$ to small cardinality or add at least 
$\lambda$ many reals. Then, by \Lemmaof{p-7}, we have
$\forces{\poP}{\xmbox{there is no }\psi^*_n\xmbox{-absolute inaccessible cardinal}}$.

By $(\calP,\emptyset)$-\RcA, it follows that there is a ground $\uniW$ of $\uniV$ \st\
$\uniW\modelof{\xmbox{there is no }\psi^*_n\xmbox{-absolute inaccessible cardinal}}$. Again 
by \Lemmaof{p-7}, this is a contradiction.
\qedoftheorem
\begin{theorem}
  \Label{p-10} Suppose that $\lambda$ is an inaccessible cardinal, $\kappa<\lambda$ is \st\
  $V_\lambda\modelof{\kappa\xmbox{ is x-large cardinal}}$, where ``x-large cardinal'' is a 
  notion of large cardinal, for which a Laver function exists. Assume also that
  $\setof{\mu<\lambda}{\mu\mbox{ is inaccessible}}$ is stationary in $\lambda$.

  Then, 
  for each of the classes $\calP$ of \pos\ considered in \theoremof{p-3}, there are 
  $\lambda_0$ with $\lambda>\lambda_0>\kappa$, and $\poP\in\calP$ with
  $\poP\subseteq V_\kappa$ \st, for a 
  $(V_\lambda,\poP)$-generic $\genG$, we have 
  \begin{equation*}
    V_{\lambda_0}[\genG]\modelof{
    \begin{array}[t]{@{}l}
      \kappa\mbox{ is a tightly$^+$ }\calP\mbox{-Laver-generically x-large cardinal}\\
      \mbox{and }\neg({\extendedcolor\calP},\emptyset)\mbox{-\RcA}
      }.\memo{\normalsize!!!}
    \end{array}
  \end{equation*}
\end{theorem}
\prf Suppose that $\calP$ is one of the classes of \pos\ considered in \theoremof{p-3}. 
Note that then, \assertof{a} or \assertof{b} of \theoremof{p-9} holds. Let $n\in\natnums$ 
be \st\ the formula ``$\kappa$ is an x-large cardinal'' is $\Sigma_n$. By the assumption, 
there is an inaccessible cardinal $\mu$ with $\lambda>\mu>\kappa$ \st\
$V_\lambda\succ V_\mu$. Let $\lambda>\lambda_0>\mu$ be the minimal cardinal \st\
$V_\lambda\models\lambda_0\xmbox{ is \jumpto{phistar}{$\psi^*_n$-absolute} inaccessible 
  cardinal}$ {\extendedcolor --- such $\lambda_0$ exists by  
  \Lemmaof{p-8}}. \memo{\normalsize!!!}
Then we have
\begin{equation}
\Label{x-6}  V_{\lambda_0}\modelof{
  \begin{array}[t]{@{}l}
    \kappa\mbox{ is an x-large cardinal}
    }
  \end{array}
\end{equation}
In $V_{\lambda_0}$, let $\poP$ be the limit of $\kappa$-iteration with appropriate support 
as described in  
\theoremof{p-3} which forces that $\kappa$ is tightly$^+$ $\calP$-Laver generically x-large 
cardinal in the generic extension of $V_{\lambda_0}$. Let $\genG$ 
be $(V_\lambda,\poP)$-generic filter. 
Then we have
$V_{\lambda_0}{\extendedcolor[\genG]}\modelof{\kappa\xmbox{ is a tightly}^+\ 
\calP\xmbox{-Laver generically x-large cardinal}}$ 
\memo{\normalsize!!!}
and 
$V_\lambda[\genG]\succ V_\mu[\genG]$ by 
\Lemmaof{p-Lg-RcA-0-2}{\ifextended\extendedcolor,\,\assertof{1}\fi}. In particular, by  
\Lemmaof{p-7}, we have
$V_{\lambda_0}[\genG]\models\mu\xmbox{ is the largest }
\psi^*_n\xmbox{-absolute inaccessible cardinal}$. 
By \theoremof{p-9}, it follows that
$V_{\lambda_0}[\genG]\modelof{\neg(\calP,\emptyset)\mbox{-\RcA}}$. 
\qedoftheorem\qedskip

The conditions of \theoremof{p-10} are satisfied {\extendedcolor by} practically all large 
cardinal notions. \memo{\normalsize !!!} 
For example, 
under the consistency of the existence of a 2-huge cardinal, the conditions of 
\theoremof{p-10} are satisfied by x-large cardinal $=$ hyperhuge cardinal 
\memo{\normalsize!!!} 
({\extendedcolor 
  see \Lemmaof{p-Lg-RcA-2} below}). Thus we obtain the following:

\begin{Cor}\Label{p-10-0}
  Under the assumption of the consistency of the existence of a 2-huge cardinal, 
  the existence of a tightly$^+$ $\calP$-Laver generically hyperhuge cardinal does not 
  imply $(\calP,\emptyset)$-\RcA\ for any class $\calP$ of \pos\ as in \theoremof{p-3}. \qed
\end{Cor}

A natural strengthening of Laver-genericity does imply the full Maximality Principle (hence 
also the full Recurrence Axiom). As 
the proof of \theoremof{p-9} suggests, such property must be formulated not in a single 
formula but as an axiom scheme. 

For a natural number $n$, we call a cardinal \ifextended\begin{itemize}\item[]\fi
$\kappa$ {\It super $C^{(n)}$-hyperhuge} if   
for any $\lambda_0>\kappa$ there are $\lambda\geq\lambda_0$ with 
$V_\lambda\prec_{\Sigma_n}\uniV$, and $j$, $M\subseteq\uniV$ \st\
$\Elembed{j}{\uniV}{M}{\kappa}$, $j(\kappa)>\lambda$, $\fnsp{j(\lambda)}{M}\subseteq M$ and
$V_{j(\lambda)}\prec_{\Sigma_n}\uniV$. \ifextended\end{itemize}\fi

\ifextended\begin{itemize}\item[]\fi
$\kappa$ is {\It super $C^{(n)}$-ultrahuge} if the condition above holds with
``$\,\fnsp{j(\lambda)}{M}\subseteq M$'' replaced by ``$\,\fnsp{j(\kappa)}{M}\subseteq M$ and 
$V_{j(\lambda)}\subseteq M$''. \iftalk\hfill{\jumpto{diagram}{\small [see the chart of 
      large large cardinals]}}\fi 
\ifextended\end{itemize}\fi

If $\kappa$ is super $C^{(n)}$-hyperhuge then it is super $C^{(n)}$-ultrahuge. This can be 
shown similarly to \Lemmaof{p-0-0}. 


We shall also say that $\kappa$ is {\It super $C^{(\infty)}$-hyperhuge} ({\It super
$C^{(\infty)}$-ultrahuge}, resp.) if it is super $C^{(n)}$-hyperhuge (super
$C^{(n)}$-ultrahuge, resp.) for all natural number $n$. 

A similar kind of strengthening of the notions of large cardinals which 
we call here ``super $C^{(n)}$'' appears also in
Boney \cite{boney}. It is called in \cite{boney}
``$C^{(n)+}$'' and the notion is considered there in connection with extendibility. 

For a natural number $n$ and an iterable class $\calP$ of \pos, 
a cardinal \ifextended\begin{itemize}\item[]\fi$\kappa$ is {\It super $C^{(n)}$- $\calP$-Laver-generically ultrahuge} 
if, for 
any $\lambda_0>\kappa$ 
and for any 
$\poP\in\calP$,
there are a $\lambda\geq\lambda_0$ 
with $V_\lambda\prec_{\Sigma_n}\uniV$, a $\calP$-name $\utpoQ$ with
$\forces{\poP}{\utpoQ\in\calP}$, \st, for $(\uniV,\poP\ast\utpoQ)$-generic $\genH$, there 
are $j$, $M\subseteq \uniV[\genH]$ with  $\Elembed{j}{\uniV}{M}{\kappa}$,
$j(\kappa)>\lambda$, $\poP$, $\genH$, ${V_{j(\lambda)}}^{\uniV[\genH]}\in M$ and
${V_{j(\lambda)}}^{\uniV[\genH]}\prec_{\Sigma_n}\uniV[\genH]$. \ifextended\end{itemize}\fi

A super $C^{(n)}$- $\calP$-Laver-generically ultrahuge cardinal 
\ifextended\begin{itemize}\item[]\fi$\kappa$ is {\It tightly 
    super $C^{(n)}$- $\calP$-Laver-generically ultrahuge}
    , if {\extendedcolor additionally} 
$\cardof{\poP\ast\utpoQ}\leq j(\kappa)$ (see \footnoteof{fn-0}) \extendedcolor{holds in the 
  definition above}. 
\memo{\normalsize!!!}
\ifextended\end{itemize}\fi

{\It Super $C^{(\infty)}$- $\calP$-Laver-generically {\extendedcolor hyperhugeness}} and 
{\It tightly super $C^{(\infty)}$- $\calP$-Laver generically {\extendedcolor hyperhugeness}}
are defined similarly to super $C^{(\infty)}$-ultrahugeness.\memo{\normalsize!!!}

Note that, in general, super $C^{(\infty)}$-hyperhugeness 
and super $C^{(\infty)}$-ultrahugeness are notions not formalizable in the language 
of \ZFC\ without introducing a new constant symbol for $\kappa$ since we need infinitely 
many $\Lin$-formulas to formulate them. Exceptions are when we 
are talking about a cardinal in a set model being with one of these 
properties like in \Lemmaof{p-Lg-RcA-2} below (and in such a case ``natural number $n$'' 
actually refers to ``$n\in\omega$''), or when we are talking about a cardinal definable in
$\uniV$ having these properties in an inner model like in \Corof{p-bedrock-7} or \Corof{p-bedrock-8}.\memox{???} In the latter case, 
the situation is formalizable with infinitely may $\Lin$-sentences.

In contrast, the super $C^{(\infty)}$-$\calP$-Laver generically ultrahugeness 
of $\kappa$ is expressible in infinitely many $\Lin$-sentences. 
This is because a $\calP$-Laver generic large 
cardinal $\kappa$ for relevant classes $\calP$ of \pos\ 
is uniquely determined as $\kappa_{\refl}$ or $2^{\aleph_0}$ (see e.g.\ \theoremof{p-1-0} 
and \theoremof{p-2}). 

A modification of the proof of \theoremof{p-Lg-RcA-0} shows the following: 

\begin{theorem}{\rm (Fuchino and Usuba \cite{laver-gen-maximum})}
  \Label{p-Lg-RcA-5} Suppose that $\calP$ is an iterable class of \pos\ and $\kappa$ 
  is {\darkred tightly} 
  \memo{correction needed if possible !!!}
  super $C^{(\infty)}$-$\calP$-Laver-generically ultrahuge. Then 
  $(\calP,\calH(\kappa))$-\RcAp\ holds. \ifextended\else\qed\fi
\end{theorem}
{\ifextended\extendedcolor
\prf Suppose that $\kappa$ is tightly super $C^{(\infty)}$-$\calP$-Laver-gen.\ ultrahuge, 
$\poP\in\calP$, and $\forces{\poP}{\varphi(\overline{a}\checked)}$ for 
an $\Lin$-formula $\varphi$ and $\overline{a}\in\calH(\kappa)$. We want to show that
$\varphi(\overline{a})$ holds in some $\calP$-ground of $\uniV$.

Let $n$ be a sufficiently large natural number \st\ the following arguments go through. In 
particular, we assume that ${V_\alpha}^\uniV\prec_{\Sigma_n}\uniV$ implies that
``$\varphi(\overline{x})$'' and ``$\forces{\cdot}{\varphi(\overline{x}\checked)}$ are 
absolute between ${V_\alpha}^\uniV$ and $\uniV$, 
and ${V_\alpha}^\uniV\prec_{\Sigma_n}\uniV$ also implies that a sufficiently large fragment 
of \ZFC\ holds in $V_\alpha$. 

Let $\utpoQ$ be a $\poP$-name \st\ $\forces{\poP}{\utpoQ\in\calP}$ and, for
$(\uniV,\poP\ast\utpoQ)$-generic $\genH$, there are a $\lambda>\kappa$ 
with 
\begin{xitemize}
\xitemA[x-Lg-RcA-9] 
  $V_\lambda\prec_{\Sigma_n}\uniV$, 
\end{xitemize}
and $j$, $M\subseteq \uniV[\genH]$ \st\ $\Elembed{j}{\uniV}{M}{\kappa}$,
$j(\kappa)>\lambda$, $\poP$, $\genH$, ${V_{j(\lambda)}}^{\uniV[\genH]}\in M$,
$\cardof{\poP\ast\utpoQ\leq j(\kappa)}$  ($<j(\kappa)$, and 
${V_{j(\lambda)}}^{\uniV[\genH]}\prec_{\Sigma_n}\uniV[\genH]$. 

Replacing $\poP\ast\utpoQ$ by an appropriate isomorphic \po\ (and replacing $\calH$ by 
corresponding filter), we may assume thet $\poP\ast\utpoQ\in {V_{j(\lambda)}}^\uniV$. 

By the choice of $n$, we 
have $V_\lambda\models\forces{\poP}{\varphi(\overline{a}\checked)}$.
$j({V_\lambda}^\uniV)={V_{j(\lambda)}}^M\prec_{\Sigma_n}M$ by elementarity of $j$,   
and 
\begin{xitemize}
\xitemA[x-Lg-RcA-9-0] 
  ${V_{j(\lambda)}}^M={V_{j(\lambda)}}^{\uniV[\genH]}$
\end{xitemize}
by the closedness of $M$. 
Since $V_\lambda\prec_{\Sigma_n}\uniV$, we have
$V_\lambda[\genH]\prec_{\Sigma_{n_0}}\uniV[\genH]$ for a still large 
enough $n_0\leq n$ by \Lemmaof{p-Lg-RcA-0-2}{\ifextended\extendedcolor,\,\assertof{1}\fi}. 
Since 
${V_{j(\lambda)}}^{\uniV[\genH]}\prec_{\Sigma_n}\uniV[\genH]$, it  
follows 
that ${V_\lambda}^{\uniV[\genH]}\prec_{\Sigma_{n_0}}{V_{j(\lambda)}}^{\uniV[\genH]}$. 
Thus 
\begin{xitemize}
\xitemA[x-Lg-RcA-10] 
  ${V_\lambda}^\uniV\prec_{\Sigma_{n_1}}{V_{j(\lambda)}}^\uniV$
\end{xitemize}
for a still large enough
$n_1\leq n_0$ by \Lemmaof{p-Lg-RcA-0-2}{\ifextended\extendedcolor,\,\assertof{2}\fi}. 

In particular, we have
${V_{j(\lambda)}}^\uniV\models\forces{\poP}{\varphi(\overline{a}\checked)}$,  and hence
$V_{j(\lambda)}[\genG]\models\varphi(\overline{a})$ where $\genG$ is the $\poP$-part of
$\genH$. Note that by \xitemAof{x-Lg-RcA-9} and \xitemAof{x-Lg-RcA-10}, $V_{j(\lambda)}$ 
satisfies a sufficiently large fragment of \ZFC. 

Thus we have
$V_{j(\lambda)}[\genH]\modelof{\,\mbox{there is a }
  \calP\mbox{-ground satisfying }\varphi(\overline{a})}$, and hence 
\begin{itemize}
\item[] ${V_{j(\lambda)}}^{\uniV[\genH]}\modelof{\,\mbox{there is a }
  \calP\mbox{-ground satisfying }\varphi(\overline{a})}$
\end{itemize}
by \Lemmaof{p-Lg-RcA-0-0}. By \xitemAof{x-Lg-RcA-9-0} and elementarity, it follows that
\begin{itemize}
\item[] ${V_{\lambda}}\modelof{\,\mbox{there is a }
  \calP\mbox{-ground satisfying }\varphi(\overline{a})}$.
\end{itemize}
Finally, this implies ${\uniV}\modelof{\,\mbox{there is a }
  \calP\mbox{-ground satisfying }\varphi(\overline{a})}$ by \xitemAof{x-Lg-RcA-9}.\\
\qedoftheorem\qedskip\fi}

The following Lemma can be proved similarly to Theorem 5c in Barbanel-DiPrisco-Tan 
\cite{barbanel-etal} (see also Theorem 24.13 in Kanamori \cite{higher-inf}). 

\begin{Lemma}{\rm (Fuchino and Usuba \cite{laver-gen-maximum})}
  \Label{p-Lg-RcA-2} If $\kappa$ is $2$-huge with the $2$-huge elementary embedding $j$, 
  that is,\ there is $M\subseteq\uniV$ \st\
  $\Elembed{j}{\uniV}{M\subseteq \uniV}{\kappa}$,  and  
  \begin{equation}
  \Label{x-Lg-RcA-5} 
    \fnsp{j^2(\kappa)}{M}\subseteq M, 
  \end{equation}
  then $V_{j(\kappa)}\modelof{\kappa
    \mbox{ is super }C^{(\infty)}\mbox{-hyperhuge cardinal\/}}$, 
  and for each $n\in\omega$,\\
  $V_{j(\kappa)}\modelof{\mbox{there are stationarily
      many super }C^{(n)}\mbox{-hyperhuge cardinals\/}}$. \qed
\end{Lemma}

The proof of the existence of Laver-function for a supercompact cardinal can be modified to 
show that super $C^{(\infty)}$-hyperhuge cardinal in $V_\mu$ has a Laver function for super
$C^{(\infty)}$-hyperhugeness (\cite{laver-gen-maximum}). Similarly to 
\theoremof{p-3} we  
obtain the following: 

\begin{theorem}{\rm (Fuchino and Usuba \cite{laver-gen-maximum})}
  \Label{p-Lg-RcA-4} \wassertof{$A$} Suppose that $\mu$ is inaccessible and 
  $\kappa<\mu$ is super $C^{(\infty)}$-hyperhuge in $V_\mu$. Let 
  $\poP=\Col(\aleph_1,\kappa)$. Then,  
  in $V_\mu[\genG]$, for any ${V_\mu,\poP}$-generic $\genG$, $\aleph_2^{V_\mu[\genG]}$
  ($=\kappa$) is tightly super $C^{(\infty)}$-$\sigma$-closed-Laver-generically hyperhuge
    and $\CH$ holds.\smallskip

  \wassert{$B$}
  Suppose that $\mu$ is inaccessible and $\kappa<\mu$ is 
  super $C^{(\infty)}$-hyperhuge with a 
  Laver function $\mapping{f}{\kappa}{V_\kappa}$ for super $C^{(\infty)}$-hyperhugeness in
  $V_\mu$.   
  If $\poP$ is the CS-iteration of length $\kappa$ for forcing \PFA\ along with $f$, then,  
  in $V_\mu[\genG]$ for any $(V_\mu,\poP)$-generic $\genG$, $\aleph_2^{V_\mu[\genG]}$
  ($=\kappa$)  
  is tightly super $C^{(\infty)}$-proper-Laver-generically hyperhuge and 
  $\continuum=\aleph_2$ holds.\smallskip

  \wassert{$B$\rlap{$'$}\,}
  Suppose that $\mu$ is inaccessible and $\kappa<\mu$ is super $C^{(\infty)}$-hyperhuge with a
  Laver function $\mapping{f}{\kappa}{V_\kappa}$ for super $C^{(\infty)}$-hyperhugeness in
  $V_\mu$.  
  If $\poP$ is the RCS-iteration of length $\kappa$ 
  for forcing \MM\ along with $f$, then,   
  in $V_\mu[\genG]$ for any $(V_\mu,\poP)$-generic $\genG$, $\aleph_2^{V_\mu[\genG]}$
  ($=\kappa$)  
  is tightly super $C^{(\infty)}$-stationary\_preserving-Laver-generically hyperhuge and 
  $\continuum=\aleph_2$ holds.\smallskip\memox{\assertof{3}, should be 
    also rewritten to assertions in the form of $V_\mu\models$ ... }

  \wassert{$\Gamma$}
  Suppose that $\mu$ is inaccessible and $\kappa$ is super $C^{(\infty)}$-hyperhuge with a 
  Laver function $\mapping{f}{\kappa}{V_\kappa}$ for super $C^{(\infty)}$-hyperhugeness in
  $V_\mu$. If 
  $\poP$ is a FS-iteration of length $\kappa$ for  
  forcing \MA\  
  along with $f$,   
  then, in $V_\mu[\genG]$ for any $(V_\mu,\poP)$-generic $\genG$, $\continuum$ ($=\kappa$) is 
  tightly super $C^{(\infty)}$-c.c.c.-Laver-generically hyperhuge, and $\continuum$ is 
  very large in $V_\mu[\genG]$.

  \wassert{$\Delta$}
  Suppose that $\mu$ is inaccessible and $\kappa$ is super $C^{(\infty)}$-hyperhuge with a 
  Laver function $\mapping{f}{\kappa}{V_\kappa}$ for super $C^{(\infty)}$-hyperhugeness in
  $V_\mu$. If 
  $\poP$ is a FS-iteration of length $\kappa$ 
  along with $f$ enumerating ``all'' \pos,   
  then, in $V_\mu[\genG]$ for 
  any $(V_\mu,\poP)$-generic $\genG$, $\continuum$ ($=\aleph_1$) is  
  tightly super $C^{(\infty)}$-all \pos-Laver-generically hyperhuge, and \CH\ 
  holds.\ifextended\else\qed\fi 
\end{theorem}
{\ifextended\extendedcolor
\memo{{\tt recurrence-axioms.tex p-Lg-RcA-4}}
\prf The proof can be done similarly to that of Theorem 5.2 in \cite{II}. In the following 
we shall only check the case \assertof{$\Delta$}. 

Suppose that $\mapping{f}{\kappa}{V_\kappa}$ is a super $C^{(\infty)}$-hyperhuge Laver 
function. 

Let $\seqof{\poP_\alpha,\utpoQ_\beta}{\alpha\leq\kappa,\beta<\kappa}$ be a FS-iteration 
defined by
\begin{xitemize}
\xitemA[p-Lg-RcA-4-a-0] $\utpoQ_\beta={}
  \left\{\,\begin{array}{@{}ll}
    f(\alpha), &\mbox{ if }f(\alpha)\mbox{ is a }\poP_\beta\mbox{-name of a \po};\\[\jot]
    \ssetof{\bbone}, &\mbox{otherwise}
  \end{array}\right.
$
\end{xitemize}
for $\beta<\kappa$.

Let $\genG$ be a $(V_\mu,\poP_\kappa)$-generic filter. Clearly
$V_\mu[\genG]\modelof{2^{\aleph_0}=\kappa=\aleph_1}$. We show that $\kappa$ is tightly 
super $C^{(\infty)}$-all \pos-Laver-generically ultrahuge in $V_\mu[\genG]$.

Suppose that $\poP$ is a \po\ in $V_\mu[\genG]$, $\kappa<\lambda_0$ and
$n\in\omega$. Let $n'>n$ be sufficiently large and let $\utpoP$ be a $\poP_\kappa$-name of
$\poP$. 

Working in $V_\mu$, we can find 
\begin{xitemize}
\xitemA[x-Lg-RcA-6-1-a] 
  $\cardof{\poP}<\lambda<\kappa^*<\lambda^*$ and $j$, $M\subseteq\uniV$
\end{xitemize}
\st\
\begin{xitemize}
\xitemA[x-Lg-RcA-6-1-0] 
  $\Elembed{j}{\uniV}{M}{\kappa}$,
\xitemA[x-Lg-RcA-6-1-1] $j(\kappa)=\kappa^*$, $j(\lambda)=\lambda^*$,
\xitemA[x-Lg-RcA-6-1-2] $\fnsp{\lambda^*}{M}\subseteq M$,
\xitemA[x-Lg-RcA-6-1-3] $V_\lambda\prec_{\Sigma_{n'}}\uniV$,
  $V_{\lambda^*}\prec_{\Sigma_{n'}}\uniV$, and 
\xitemA[x-Lg-RcA-6-1-4] $j(f)(\kappa)=\utpoP$ 
\end{xitemize}
by definition of $f$. 


By elementarity (and by the definition \xitemAof{p-Lg-RcA-4-a-0} of $\poP_\kappa$),
\begin{xitemize}
\xitemA[x-Lg-RcA-6-2] 
  $j(\poP_\kappa)\ \sim_{\poP_\kappa}\ (\poP_\kappa\ast\utpoP)\ast\utpoR$ 
\end{xitemize}
for a $(\poP_\kappa\ast\utpoP)$-name $\utpoR$ of a \po.
Note that $(\poP_\kappa\ast\utpoP\ast\utpoR)/\genG$ corresponds to a \po\ of the form
$\poP\ast\utpoQ$.

Let $\genH^*$ be $(\uniV,(\poP_\kappa\ast\utpoP)\ast\utpoR)$-generic filter 
with $\genG\subseteq\genH^*$. $\genH^*$ corresponds to a $(\uniV,j(\poP_\kappa))$-generic 
filter $\genH\supseteq\genG$ via the equivalence \xitemAof{x-Lg-RcA-6-2}.

Let $\tilde{j}$ be defined by
\begin{xitemize}
\xitemA[x-Lg-RcA-6-3] 
  $\mapping{\tilde{j}}{\uniV[\genG]}{M[\genH]}$;\quad $\uta^\genG\mapsto j(\uta)^\genH$
\end{xitemize}
for all $\poP_\kappa$-names $\uta$. 

A standard proof shows that $f$ is well-defined, and
$\Elembed{j}{\uniV[\genG]}{M[\genH]}{\kappa}$.   
By \xitemAof{x-Lg-RcA-6-1-1} and \xitemAof{x-Lg-RcA-6-1-2}, we have
$\tilde{j}\imageof \tilde{j}(\lambda)=j\imageof j(\lambda)\in M[\genH]$.  
Since $\genH\in M[\genH]$, the $(\uniV[\genG],\poP\ast\utpoQ)$-generic filter corresponding 
to $\genH$ is also in $M[\genH]$. 

By \xitemAof{x-Lg-RcA-6-1-a}, \xitemAof{x-Lg-RcA-6-1-3}, by the choice of $n'$,
and by 
\Lemmaof{p-Lg-RcA-0-2}{\ifextended\extendedcolor,\,\assertof{1}\fi}, we have 
${V_\lambda}^{\uniV[\genG]}\prec_{\Sigma_n}\uniV[\genG]$\quad and\quad
${V_{\tilde{j}(\lambda)}}^{\uniV[\genH]}
={V_{\lambda^*}}^{\uniV[\genH]}\prec_{\Sigma_n}\uniV[\genH]$.

Since $\poP$ and $n$ were arbitrary, this shows that $\kappa$ is tightly 
super $C^{(\infty)}$-all \pos-Laver-generically ultrahuge in $V_\mu[\genG]$.
\qedoftheorem\qedskip\fi}

\section{Toward the Laver-generic Maximum}
\Label{sec:4} Besides \theoremof{p-Lg-RcA-0} and \theoremof{p-Lg-RcA-5}, 
\memox{{\normalsize!!!}The  
  theorem about ``super $C^{\infty}$ ...'' implies RcA+}, 
we also have some other advantages of assuming the existence of $\calP$-Laver-generically 
hyperhuge cardinal or even its ``super $C^{(\infty)}$'' version. One of them is that 
they imply the existence of the (set-generic) bedrock (see below for definition); another is 
that we know the exact  
consistency strength of these principles. 

For a class $\calP$ of \pos, 
a cardinal $\kappa$ is {\It tightly $\calP$-generic hyperhuge} if for any 
$\lambda>\kappa$, there is $\poQ\in\calP$ \st\ for a $(\uniV,\poQ)$-generic $\genH$, there 
are $j$, $M\subseteq\uniV[\genH]$ \st\ $\Elembed{j}{\uniV}{M}{\kappa}$,
$\lambda<j(\kappa)$, 
$\cardof{\poQ}\leq j(\kappa)$, and
$j\imageof{j(\lambda)}, \genH\in M$. 


Note that, for any class $\calP$ of \pos\ with $\ssetof{\bbone}\in\calP$, 
the hyperhugeness of a cardinal $\kappa$ implies its tightly $\calP$-generically 
hyperhugeness. 
Likewise, if $\calP$ is iterable then, the tightly $\calP$-Laver-generic hyperhugeness of 
$\kappa$ implies its tightly $\calP$-generically 
hyperhugeness.\nopagebreak
\medskip
\nopagebreak\linebreak
\noindent
\mbox{}\hfill
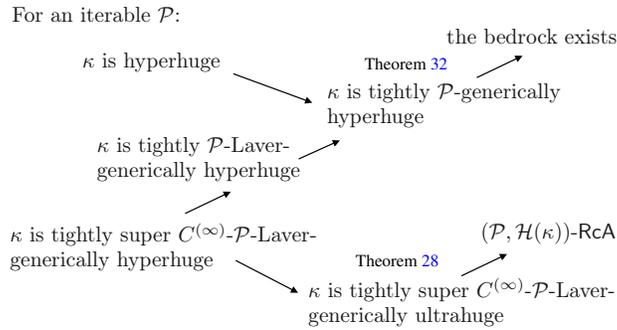
\includegraphics[width=0.7\textwidth]{reflection\_and\_recurrence-Janos-Festschrift-fig1.pdf}\hfill\mbox{}
\nopagebreak
\vspace{-24.6ex}
\nopagebreak\linebreak
\mbox{}\hspace{0.55\textwidth}{\scalebox{0.76}{\footnotesize \theoremof{p-12}}}
\vspace{13.9ex}
\nopagebreak\linebreak
\mbox{}\hspace{0.54\textwidth}{\scalebox{0.76}{\footnotesize \theoremof{p-Lg-RcA-5}}}
\vspace{5ex}
\medskip

Usuba \cite{usuba} proved that the grounds of $\uniV$ are downward directed (\wrt\ subclass 
relation) for class many grounds (this is formalizable by virtue of \theoremof{p-intro-1}). 
More concretely
\begin{theorem}{\rm (Theorem 1.3 in Usuba \cite{usuba})} 
  \Label{p-11}
  For any collection of grounds of
  $\uniV$, indexed by a set of parameters (in the sense of \theoremof{p-intro-1}), there is 
  a ground which is included in all grounds in the collection. \qed
\end{theorem}
From this theorem, it follows that the {\It mantle}, i.e., the intersection
of all grounds is 
a model of \ZFC. In \cite{usuba}, it is proved that the mantle is a ground and hence it is 
the {\It bedrock}, i.e., the smallest ground of $\uniV$ provided that there exists a 
hyperhuge cardinal(Theorem 1.6 in \cite{usuba}). Later the assumption of the existence of a 
hyperhuge 
cardinal in this theorem is weakened to the existence of an extendible cardinal (Theorem 
1.3 in Usuba \cite{usuba2}).  

In \cite{laver-gen-maximum}, we obtained the following generalization of Theorem 1.6 in 
\cite{usuba}. In the following theorems, tightness of a $\calP$-generic Large cardinal is 
defined similarly to the tightness of $\calP$-Laver generic large cardinal. Note that 
the (in many cases unique (\theoremof{p-2})) 
tightly $\calP$-Laver generic large cardinal as well as 
corresponding genuine large 
cardinals are tightly $\calP$-generic large cardinals by 
definition. 

\begin{theorem}{\rm(Fuchino and Usuba \cite{laver-gen-maximum})}\Label{p-12} If there is a 
  tightly $\calP$-generically hyperhuge cardinal $\kappa$, then the mantle is a ground of
  $\uniV$. In particular it is the bedrock.
\end{theorem}
\noindent
{\it A sketch of the proof.\/} The overall structure of the structure of the proof is just 
the same as that of Theorem 1.6 in \cite{usuba} or Theorem 1.3 in \cite{usuba2}.

We call a ground $\uniW$ of $\uniV$ a $\LE\kappa$-ground if there is $\poP\in W$ with
$\cardof{\poP}^\uniV\leq\kappa$ and a $(\uniW,\poP)$-generic $\genG$ \st\
$\uniV=\uniW[\genG]$. Let
\begin{equation}\Label{x-7}
  \overline{W}=\bigcap\setof{W}{W\mbox{ is a }\LE\kappa\mbox{-ground}}.
\end{equation}
By \theoremof{p-11}, there is a ground $\uniW\subseteq\overline{\uniW}$.
For such $\uniW$ it is enough to show that actually $\overline{\uniW}\subseteq\uniW$ holds.

Let $\poS\in\uniW$ be a \po\ with cardinality $\mu$ (in $\uniV$) \st\ there is a
$(\uniW,\poS)$-generic $\genF\in\uniV$  with $\uniV=\uniW[\genF]$. \Wolog, $\mu\geq\kappa$.

By \theoremof{p-intro-1}, there is $r\in\uniV$ \st\ $\uniW=\Phi(\cdot,r)^\uniV$. 

Let $\theta\geq\mu$ be \st\ $r\in V_\theta$, and for a sufficiently large natural 
number $n$, we have 
${V_\theta}^\uniV\prec_{\Sigma_n}\uniV$.
By the choice of $\theta$, 
$\Phi(\cdot,r)^{{V_\theta}^\uniV}=\Phi(\cdot,r)^\uniV\cap {V_\theta}^\uniV
  =\uniW\cap {V_\theta}^\uniV={V_\theta}^\uniW$. Let 
$\poQ\in\calP$ \st\ for $(\uniV,\poQ)$-generic $\genH$, there are $j$,
$M\subseteq\uniV[\genH]$ with 
$\Elembed{j}{\uniV}{M}{\kappa}$, $\theta<j(\kappa)$,
$\cardof{\poQ}\leq j(\kappa)$,
${V_{j(\theta)}}^{\uniV[\genH]}\subseteq M$, and 
$\genH$, $j\imageof{j(\theta)}\in M$. 

Using this $j$ we can show that 
${V_\theta}^{\overline{\uniW}}\subseteq{V_\theta}^{\uniW}$ holds (this part of the proof is 
quite involved, for the details, the reader is referred to \cite{laver-gen-maximum}). 
Since $\theta$ can be arbitrary large, It follows that $\overline{\uniW}\subseteq\uniW$.
\qedoftheorem\qedskip

Analyzing the details of the proof of \theoremof{p-12} which we omitted from our present 
exposition, we also obtain the following result with many surprising consequences:

\begin{theorem}{\rm(Fuchino and Usuba \cite{laver-gen-maximum})}
  \Label{p-bedrock-2} Suppose that $\calP$ is any class of \pos.
  If $\kappa$ is a tightly $\calP$-generically hyperhuge cardinal, then $\kappa$ is a 
  hyperhuge cardinal in the bedrock $\overline{W}$ of $\uniV$. \qed
\end{theorem}

The following equiconsistency results are immediate consequences of the theorem above:
\begin{Cor}
  \Label{p-bedrock-5} Suppose that $\calP$ is the class of all \pos. 
  Then the following theories are 
  equiconsistent:\smallskip
  
  \wassert{a} \ZFC\ $+$ ``there is a hyperhuge cardinal''.

  \wassert{b} \ZFC\ $+$ ``there is a tightly $\calP$-Laver generically hyperhuge cardinal''.

  \wassert{c} \ZFC\ $+$ ``there is a tightly $\calP$-generically hyperhuge cardinal''.

  \wassert{d} \ZFC\ $+$ ``the bedrock\/ $\overline{\uniW}$ exists and $\omega_1$ is a 
  hyperhuge cardinal in $\overline{\uniW}$''.\qed 
\end{Cor}

\begin{Cor}
  \Label{p-bedrock-6} Suppose that $\calP$ is one of the following classes of \pos: 
  all semi-proper \pos; all proper \pos; all ccc \pos; all $\sigma$-closed \pos. 
  Then the following theories are 
  equiconsistent:\smallskip
  
  \wassert{a} \ZFC\ $+$ ``there is a hyperhuge cardinal''.

  \wassert{b} \ZFC\ $+$ ``there is a tightly $\calP$-Laver generically hyperhuge cardinal''.

  \wassert{c} \ZFC\ $+$ ``there is a tightly $\calP$-generically hyperhuge cardinal''.

  \wassert{d} \ZFC\ $+$ ``the bedrock\/ $\overline{\uniW}$ exists and $\kappa_\refl$ is a 
  hyperhuge cardinal in $\overline{\uniW}$''. \qed
\end{Cor}

These equiconsistency results are quite remarkable when we remember that equiconsistency of 
axioms like \PFA, \MM, \MM$^{++}$ etc.\ are unknown at the moment.

The tightness of generic large cardinals was originally 
thought as a technical condition when it was introduced to overcome the difficulty in the 
proof of \theoremof{p-1-0}, \assertof{2}. 
The (proofs of) theorems \ref{p-12} and \ref{p-bedrock-2} and their consequences, some of 
which we are presenting here as their corollaries, suggest that this notion is much more 
intrinsic than merely a technicality.

A slight modification of the proofs of the theorems above also show the following. Note that 
as we already noticed, super-$C^{(\infty)}$-large cardinal is not formalizable in the 
language of \ZFC. However, the assertions \assertof{a} and \assertof{b} in the following 
\Corof{p-bedrock-7} and \Corof{p-bedrock-7} can be formulated as schemes of sentences in
$\Lin$.  

\begin{Cor}
  \Label{p-bedrock-7} Suppose that $\calP$ is the class of all \pos. 
  Then the following theories are 
  equiconsistent:\smallskip
  

  \wassert{a} \ZFC\ $+$ ``$c$ is a super $C^{(\infty)}$ hyperhuge cardinal'' where $c$ is a 
  new constant symbol but ``... is super $C^{(\infty)}$ hyperhuge ...'' is formulated in an 
  infinite collection of formulas in $\Lin$. \smallskip

  \wassert{b} \ZFC\ $+$ ``there is a tightly super $C^{(\infty)}$-$\calP$-Laver generically 
  hyperhuge cardinal''. \smallskip

  \wassert{c} \ZFC\ $+$ ``the bedrock\/ $\overline{\uniW}$ exists and $\omega_1^\uniV$ is a 
  super $C^{(\infty)}$-hyperhuge cardinal in $\overline{\uniW}$''. \qed 
\end{Cor}

\begin{Cor}
  \Label{p-bedrock-8} Suppose that $\calP$ is one of the following classes of \pos: 
  all semi-proper \pos; all proper \pos; all ccc \pos; all $\sigma$-closed \pos. 
  Then the following theories are 
  equiconsistent:\smallskip
  
  \wassert{a} \ZFC\ $+$ ``$c$ is a super $C^{(\infty)}$ hyperhuge cardinal'' where $c$ is a 
  new constant symbol but ``... is super $C^{(\infty)}$ hyperhuge ...'' is formulated in an 
  infinite collection of formulas in $\Lin$. \smallskip

  \wassert{b} \ZFC\ $+$ ``there is a tightly super $C^{(\infty)}$-$\calP$-Laver generically 
  hyperhuge cardinal''. \smallskip

  \wassert{c} \ZFC\ $+$ ``the bedrock\/ $\overline{\uniW}$ exists and ${\kappa_\refl}^\uniV$ is a super
  $C^{(\infty)}$-hyperhuge cardinal in $\overline{\uniW}$''. \qed
\end{Cor}

Finally, we move to the promised proof of \theoremof{p-1-0},\,\assertof{3}.

\begin{Cor}
  \Label{p-bedrock-3} Suppose that $\calP$ is an arbitrary class of \pos\ and 
  $\kappa$ is a tightly $\calP$-generically hyperhuge cardinal.  Then\smallskip

  \wassert{1} there are 
  cofinally many huge \rlap{cardinals.}\smallskip

  \wassert{2} \SCH\ holds above some cardinal. 
\end{Cor}
\prf Suppose that $\kappa$ is a tightly $\calP$-generically hyperhuge cardinal. 
By \theoremof{p-12}, there is the bedrock $\overline{\uniW}$ and $\kappa$ is hyperhuge 
cardinal in $\overline{\uniW}$. \smallskip

\assertof{1}: 
Since the existence of a hyperhuge cardinal implies the 
existence of cofinally many huge cardinals (it is easy to show that the target $j(\kappa)$ 
of hyperhuge embedding for a sufficiently large inaccessible $\lambda$ is a huge cardinal), 
there are cofinally many huge cardinals in 
$\overline{\uniW}$. Since $\uniV$ is attained by a set forcing starting 
from $\overline{\uniW}$, a final segment of these huge cardinals survive in $\uniV$. 
\smallskip

\assertof{2}: By Theorem 20.8 in \cite{millennium-book}, \SCH\ holds above $\kappa$ in
$\overline{\uniW}$. Since $\uniV$ is a set generic extension of $\overline{\uniW}$. \SCH\ 
should hold above some cardinal $\mu\geq\kappa$. 
\qedofCor\qedskip

For iterable stationary preserving $\calP$ containing all proper \pos, 
\Corof{p-bedrock-3},\,\assertof{2} 
holds already under the $\calP$-Laver-generic supercompactness of $\kappa$. The reason is 
that in such case \PFA\ holds by \theoremof{p-1}, and by Viale \cite{viale}, \SCH\ follows 
from it. \qedskip

\noindent
{\it Proof of {\rm\theoremof{p-1-0}},\,\assertof{3}.}\quad
Let $\lambda$ and $\poQ$ be \st\ 
\begin{xitemize}
\xitem[x-Lg-RcA-8-a-0] 
  $\lambda>2^{\aleph_0}$,
  $\kappa$ and $\lambda$ is large enough \st\ \SCH\ holds above some $\mu<\lambda$ (this is 
  possible by \Corof{p-bedrock-3} ,\,\assertof{2}, and it is here that we need a strong 
  property like the Laver generic hyperhugeness of $\kappa$),
\xitem[x-Lg-RcA-8-a-1] 
  $\poQ$ is positive elements of a complete Boolean algebra, 
  and, 
\xitem[x-Lg-RcA-8-a-2] 
  for $(\uniV,\poQ)$-generic $\genH$, there are $j$, $M\subseteq\uniV[\genH]$ \st\
  \wassertof{1} $\Elembed{j}{\uniV}{M}{\kappa}$,
  \wassertof{2} $j(\kappa)>\lambda$, \wassertof{3} $\cardof{\poQ}\leq j(\kappa)$,
  and \wassertof{4} ${V_{j(\lambda)}}^{\uniV[\genH]}\subseteq M$. 
\end{xitemize}

By \xitemof{x-Lg-RcA-8-a-1}, 
each $\poQ$-name $\utilde{r}$ of a 
real 
corresponds to a mapping $\mapping{f}{\omega}{\poQ}$. 
By \xitemof{x-Lg-RcA-8-a-0} and by \xitemof{x-Lg-RcA-8-a-2}, 
\assertof{3}, 
there are at most $j(\kappa)$ many such mappings. Thus 
we have $\uniV[\genH]\modelof{\continuum\leq j(\kappa)}$, By \xitemof{x-Lg-RcA-8-a-2}, 
\assertof{4}, it follows 
$M\modelof{\continuum\leq j(\kappa)}$. By elementarity, it follows that
$\uniV\modelof{\continuum\leq\kappa}$.
\qedof{\theoremof{p-1-0},\,\assertof{3}}\qedskip

Returning to \assertof{$E$} and \assertof{$Z$} at the end of \sectionof{sec:2}, we now know 
that the existence of a super $C^{(\infty)}$-stationary preserving-Laver-generically  
hyperhuge cardinal implies \assertof{$E$} (actually it even implies
$(\mbox{stationary preserving \pos, }\calH(\kappa_\refl))$-\RcAp), and
that the existence of a super $C^{(\infty)}$-all \pos-Laver-generically  
hyperhuge cardinals implies \assertof{$Z$} (actually it even implies
$(\mbox{all \pos, }\calH(2^{\aleph_0}))$-\RcAp). These two scenarios are not compatible 
since the former implies $\continuum=\aleph_2$ while the latter implies \CH.

However, with the following axiom, \assertof{$E$} is reconciled with a fragment of 
\assertof{$Z$}: 
\begin{itemize}
\item[{\darkred \LGM:}]\qquad\quad the continuum is tightly 
  super $C^{(\infty)}$-stationary preserving-Laver generically hyperhuge and 
  there is a ground $\uniW$ of $\uniV$ \st\ the continuum is tightly 
  super $C^{(\infty)}$-all \pos-Laver generically hyperhuge in $\uniW$. 
\end{itemize}

This combination of Laver-genericity and ``ground''\footnote{I am using the word ``ground'' 
  here as as an adjective contrasting with the word ``generic'' in 
  ``generic large cardinal''.} Laver-genericity above implies \assertof{$Z^+$} mentioned at 
the end of \sectionof{sec:2} (by \theoremof{p-Lg-RcA-5}). As \assertof{$Z^+$} represents in 
a sense  the maximal amount of available strengthening of recurrence, I would like to 
choose the name Laver-Generic Maximum (\LGM) for it. 

If we admit that Recurrence Axioms, Maximal Principles and Laver-genericity are natural 
requirements, we should be also ready to accept \ZFC\ $+$ \LGM\ as a natural candidate of 
the extension of \ZFC.  

By \theoremof{p-1}, \LGM\ implies the double plus version of Martin's Maximum (\MM$^{++}$) 
and hence all the consequences of it including 
$\continuum=\aleph_2$.

\phantomsection
\Label{CichonM}
By \theoremof{p-12}, \LGM\ implies that there is the bedrock. 
So by \theoremof{p-Lg-RcA-5}, \LGM\ implies 
\assertof{$Z^+$} on page \pageref{z+}. \assertof{$Z^+$} implies that if some statement
$\varphi$ is forcable by a stationary preserving \po, then for any $A\in\calH(\aleph_2)$, 
there is a semi-proper-ground $\uniW$ of $\uniV$ \st\ $A\in\uniW$ 
and $\uniW\models\varphi$. In particular, Cichoń's Maximum \cite{GKS}, \cite{GKMS} is a 
phenomena in many semi-proper-grounds in this sense. 
Note that, by \Corof{p-bedrock-3},\,\assertof{1} the forcing argument for $\varphi$ may 
even utilize class 
many huge cardinals (e.g.\ the proof in \cite{GKS} uses four strongly compact cardinals).
{\footnote{\extendedcolor\memo{\normalsize!!!}In \cite{laver-gen-maximum}, we even show that 
  \LGM\ implies that there are stationarily many super-$C^{(\infty)}$-hyperhuge cardinals.}}

Even in the case that the forcing to show the consistency of $\varphi$  
is not stationary preserving, we can still find some ground $\uniW$ of $\uniV$ which 
satisfies $\varphi$. \memox{{\normalsize !!!}existence of the bedrock}

\phantomsection
\Label{LGM}
Thus \ZFC\ $+$ \LGM\ is a very strong axiom system which integrates practically all 
statements into itself as far as these statements can be proved to be consistent by way of 
forcing and/or methods of inner models or some combination of them. Against this backdrop, 
we want to call the axiom system \LGM\ (or possibly some 
further extension of it in the future) the {\It Laver Generic Maximum}. 

The consistency and equiconsistency of \LGM\ is easily established: we start from a model 
with two super $C^{(\infty)}$ hyperhuge cardinals $\kappa_0<\kappa_1$. We force 
$\kappa_0$ to be tightly super $C^{(\infty)}$-all \pos-Laver generically hyperhuge 
(\theoremof{p-Lg-RcA-4}, \assertof{$\Delta$}). We then force make $\kappa_1$ to be tightly 
super $C^{(\infty)}$-stationary preserving-Laver generically hyperhuge 
(\theoremof{p-Lg-RcA-4}, \assertof{$B'$}). 

By \theoremof{p-bedrock-2} the consistency strength of \LGM\ can be proved to be 
equivalent with that of two super $C^{(\infty)}$ hyperhuge cardinals (which may be 
formulated by using two new constant symbols). 

\section{More about consistency strength}
\Label{sec:5}
When I was a student of Janos Makowsky in the early 1980s at the Free University 
of Berlin,  
one of the  
papers he was preparing then was \cite{vopenka} in which the effects of Vopěnka 
Principle on properties of model theoretic logics is studied. I remember that Janos gave a 
talk on this subject in the (West) Berliner 
Logic Colloquium. Back then, I was still living in a set theory of consistency strength 
way below a measurable cardinal, and could not begin with the material of his paper at all 
because of the vertiginous consistency strength of the Vopěnka Principle.  

Janos left Berlin before I wrote up my diploma thesis on abstract elementary classes which 
was the subject Janos gave me; all the assertions I proved in the thesis remained in the 
consistency strength of \ZFC. \hfill\medskip\smallskip\Hypertarget{diagram}{}
\nopagebreak\linebreak
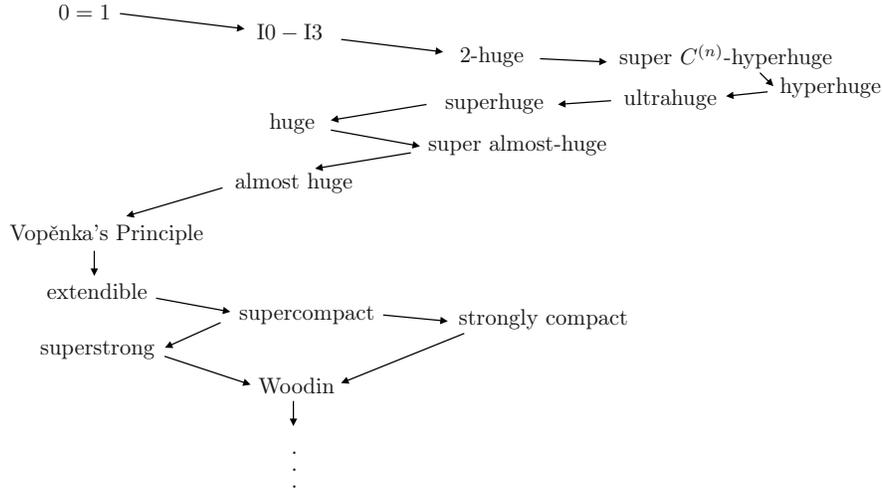
\noindent
\includegraphics[width=\textwidth]{reflection\_and\_recurrence-Janos-Festschrift-upper-half-of-higher-inf.pdf}\smallskip

The consistency strength of my set theoretic world view reached the realm of one supercompact 
cardinal when I wrote \cite{openly-gen} in the early 1990s in which some consequences of
$\MA^+(\sigma\mbox{-closed})$ were discussed.
However, it is only quite recently that I caught up with Janos definitively (at least  in terms 
of 
consistency strength) when I considered in \cite{laver-gen-maximum} the  
super $C^{(\infty)}$-$\calP$-Laver-generically hyperhuge cardinals whose 
consistency strength is (demonstrably --- see \sectionof{sec:4}) strictly between hyperhuge 
and 2-huge. 

In the meantime, active research on abstract model classes is resumed and Janos's 
\cite{vopenka} begins to attract the attention of young logicians.  
For example, Janos's \cite{vopenka} was recently cited in Boney \cite{boney} which was 
already mentioned in  
\sectionof{sec:4}\memox{\large Check!}. Boney cites the main theorem of \cite{vopenka} as 
Fact 3.12 in his paper 
and comments in allusion to  \
Aki Kanamori's comment on Kunen's inconsistency proof in \cite{higher-inf} that Vopěnka's 
Principle ``\,`rallies at least to force a veritable Götter\-dömmerung' for compactness 
cardinals for logics.'' The gap between Vopěnka's Principle and a huge cardinal 
Boney mentions in connection with this \glqq götterdämmerigen\grqq\ statement seems to 
have some 
resemblance to the discrepancy between usual Laver-generic large cardinal axioms 
and the super $C^{(\infty)}$-Laver-generic large cardinals. \memox{\large !!!} 

Now one of the urgent items in my to-do-list is to check Janos's \cite{vopenka} as well as
\cite{makowsky-XVIII}, \cite{makowsky-XX}, \cite{makowsky-mundici} 
more carefully
to find out further possible connections of his results to the context I described above. 

Is this also an instance of (the eternal?) recurrence?  

\ifextended
\phantomsection
\addcontentsline{toc}{section}{Acknowledgment} 
\fi
\begin{acknowledgement}
  A part of the material of the chapter was presented at the RIMS set theory workshop 
  2023 in Kyoto.  The author would like to thank Gunter Fuchs and Joan Bagaria for valuable 
  comments at the talk. The author also would like to thank Andrés Villaveces and the 
  anonymous referee for their 
  critical remarks on the first draft of the chapter. The draft was largely rewritten 
  according to their suggestions. The author hopes very much that the present 
  revised version of the chapter is now readable for a broader audience.

  One of the suggestions of the referee was that the terminology of Laver-genericity should 
  be simplified. We decided however that, at least for the present article, the cumbersome and 
  surely annoying expressions like ``super $C^{(\infty)}$ tightly $\calP$-Laver generically 
  hyperhuge'' will be kept as they are. The main reason for the decision is that this clumsy 
  way of saying does express subtle differences which might be crucial for the further 
  development of the theory of Laver-genericity.

  The theory of Laver-genericity is indeed still under transition. When the tightness of 
  Laver genericity was introduced in \cite{II}, the author thought that it is only a 
  technical condition. Later, in \cite{laver-gen-maximum}, we realized that the tightness 
  is {\it the} central notion in connection with the existence of bedrock (see \theoremof{p-12}). 
  Thus it is proved to be a right decision that we did not change the definition of Laver 
  genericity to make the condition of tightness a part of its definition. 

  An attempt is currently being made to integrate the theory of Laver-genericity into the 
  larger context of the set theory of Hugh Woodin (or possibly other way round). The author 
  hopes that, in a future expository article, when this project has been accomplished, the 
  whole theory can be represented in a nice and simplified choice of terminology. \medskip

  The author's research was supported by Kakenhi Grant-in-Aid 
  for Scientific Research (C) 20K03717.

\end{acknowledgement}


\newcommand{\bysame}[1]{\underline{\phantom{#1}}}%

\hfill \begin{CJK}{UTF8}{ipxm} 渕野 昌 (ふちの さかえ)
\end{CJK}
\end{document}